\documentclass[12pt,draftcls,onecolumn]{IEEEtran}
\ifCLASSINFOpdf
\else
\fi
\usepackage{amsfonts}
\usepackage{amsmath}
\usepackage{bm}
\usepackage{wrapfig}
\usepackage{graphicx}
\usepackage{epstopdf}
\usepackage{xcolor}
\usepackage[ruled,vlined,linesnumbered]{algorithm2e}
\usepackage{algpseudocode}
\usepackage{tabularx}

\newcommand{\argmax}{\textnormal{argmax}}
\newtheorem{assumption}{Assumption}

\newtheorem{theorem}{Theorem}
\newtheorem{lemma}{Lemma}
\newtheorem{proposition}{Proposition}

\newtheorem{remark}{Remark}

\begin{document}
%
\title{Optimal Control of Endo-Atmospheric Launch Vehicle Systems: Geometric and Computational Issues}
%
%
%

\author{Riccardo~Bonalli,
        Bruno~H\'{e}riss\'{e}
        and~Emmanuel~Tr\'{e}lat
\thanks{R. Bonalli is with Sorbonne Universit\'e, Universit\'e Paris-Diderot SPC, CNRS, Inria, Laboratoire Jacques-Louis Lions, \'equipe CAGE, F-75005 Paris, France and ONERA, DTIS, Universit\'e Paris Saclay, F-91123 Palaiseau, France, e-mail: rbonalli@stanford.edu, riccardo.bonalli@etu.upmc.fr, riccardo.bonalli@onera.fr.}
\thanks{B. H\'{e}riss\'{e} is with ONERA, DTIS, Universit\'e Paris Saclay, F-91123 Palaiseau, France, e-mail: bruno.herisse@onera.fr.}
\thanks{E. Tr\'{e}lat is with Sorbonne Universit\'e, Universit\'e Paris-Diderot SPC, CNRS, Inria, Laboratoire Jacques-Louis Lions, \'equipe CAGE, F-75005 Paris, France, e-mail: emmanuel.trelat@sorbonne-universite.fr.}
}

\maketitle

\begin{abstract}
In this paper we develop a geometric analysis and a numerical algorithm, based on indirect methods, to solve optimal guidance of endo-atmospheric launch vehicle systems under mixed control-state constraints. Two main difficulties are addressed. First, we tackle the presence of Euler singularities by introducing a representation of the configuration manifold in appropriate local charts. In these local coordinates, not only the problem is free from Euler singularities but also it can be recast as an optimal control problem with only pure control constraints. The second issue concerns the initialization of the shooting method. We introduce a strategy which combines indirect methods with homotopies, thus providing high accuracy. We illustrate the efficiency of our approach by numerical simulations on missile interception problems under challenging scenarios.

\end{abstract}


%
\IEEEpeerreviewmaketitle

\section{Introduction} \label{secIntro}

\subsection{Optimal Guidance of Launch Vehicle Systems}

\IEEEPARstart{G}{uidance} of autonomous launch vehicle systems towards rendezvous regions is a complex task, often considered in aerospace applications. It can be modeled as an optimal control problem with the objective of finding a control law enabling the vehicle to join some target region considering prescribed constraints as well as performance criteria. The rendezvous region may be static as well as a moving point if, for example, the mission consists of reaching a maneuvering goal. Then, an important challenge consists of developing analysis and algorithms able to provide \textit{high numerical precision} for optimal trajectories, considering rough onboard processors, i.e., reduced computational capabilities.

In the engineering community, one of the most widespread approaches to solve such kind of task relies on \textit{explicit guidance laws} (see, e.g., \cite{LinTsai,Lin,shinar1991optimal,morgan2011minimum,Nahshon}). They correct errors coming from perturbations and misreading of the system. Nonetheless, trajectories induced by guidance laws are usually not optimal because of some approximations that are required to develop a closed-form expression. On the other hand, computation of trajectories is often achieved by adopting \textit{direct methods} (see, e.g., \cite{Paris,RossOld,Ross,petersen2013model,weiss2015model}). These techniques consist of discretizing each component of the optimal control problem (the state, the control, etc.) to reduce it to a nonlinear constrained optimization problem. A high degree of robustness is provided while, in general, no deep knowledge of properties related to the structure of the dynamical system is needed, making these methods particularly easy to use in practice. However, their efficiency is proportional to the computational load which often obliges to use them offline.

Good candidates to deal with onboard processing of optimal trajectories are \textit{indirect methods} (see, e.g., \cite{calise1979singular,calise1998design,Lu1,Lu3,Pontani}). They leverage necessary conditions for optimality coming from the \textit{Pontryagin Maximum Principle} (PMP) (see, e.g., \cite{pontryagin1987mathematical,lee1967foundations}) to wrap the optimal guidance problem into a two-point boundary value problem, leading to accurate and fast algorithms (see, e.g. \cite{betts1998survey}). The advantages of indirect methods, whose more basic version is known as \textit{shooting method}, are their extremely good numerical accuracy and the fact that, when they converge, convergence is very quick. However, \textit{initializing} indirect methods is a challenging task. Moreover, further methodological difficulties arise in designing algorithms that are based on indirect methods.

\subsection{Additional Methodological Issues: Euler Coordinates Singularities coming from Mixed Control-State Constraints}

Obtaining efficient solutions for optimal guidance may oblige to consider both demanding performance criteria and possible onerous missions to accomplish. Since, in this situation, the vehicle is subject to several strong mechanical strains, some stability constraints must be imposed, which are modeled as \textit{mixed control-state constraints}. These optimal control problems are more difficult to tackle by the PMP (see, e.g., \cite{bryson1975applied,hartl1995survey,de2001maximum,clarke2010optimal}). Indeed, further Lagrange multipliers appear, for which obtaining useful information may be arduous and has been the object of many studies in the existing literature (see, e.g., \cite{jacobson1971new,maurer1977optimal,bonnard2003optimal,bonnans2007well,arutyunov2011maximum}).

A widespread approach in aeronautics to avoid to deal with these particular mixed control-state constraints consists of reformulating the original guidance problem using some \textit{local Euler coordinates}, under which the structural constraints become pure control constraints (see, e.g., \cite{bonnard2003optimal}; we discuss this change of coordinates in Section \ref{localPMP}). The transformation allows to consider the usual PMP, and then, classical shooting methods. However, Euler coordinates are not global and their singularities prevent from solving all reachable configurations, reducing the number of feasible missions.

\subsection{Statement of Contributions}

The main objective of this paper consists of designing a numerical strategy based on indirect methods to solve optimal guidance of endo-atmospheric launch vehicle systems. This strategy is able to provide global solutions lying in the configuration manifold by tackling the presence of Euler coordinates singularities introduced by mixed control-state constraints. The contribution is twofold: we first provide a geometric analysis of necessary conditions for optimality from which we derive a numerical scheme ensuring convergence of indirect methods when mixed control-state constraints are considered. The advantage of this strategy is that solutions satisfying mixed control-state constraints can be found by merely employing usual shooting methods.

Specifically, our contributions go as follows: \\

\noindent 1) \textit{Geometric analysis of necessary conditions for optimality:} The solution that we propose to bypass the problem of Euler coordinates singularities consists of reformulating the optimal guidance problem within an intrinsic viewpoint, using geometric control (it does not seem that this framework has been investigated in the optimal guidance context so far).

We build additional local coordinates that cover the singularities of the previous ones (see Section \ref{secIntro}.B) and in which the mixed control-state constraints can be expressed as pure control constraints (see Section \ref{secIntro}.B) as well. Moreover, these two sets of local coordinates form an atlas of the configuration manifold and we prove, by using geometric control techniques, that the local PMP formulations in these charts, which have only pure control constraints, are (locally) equivalent to the global PMP formulation with mixed constraints. This justifies the implementation of indirect methods to solve the original problem by employing classical shooting algorithms on the two local problems (with pure control constraints).

We stress the fact that the introduction of these particular local coordinates provides, in turn, two main benefits. On one hand, there is no limit on the feasible missions that can be simulated, and, on the other hand, the optimal guidance problem is not conditioned by multipliers depending on mixed constraints, then, standard shooting or multi-shooting methods can be easily put in practice. This is at the price of changing chart (that is, local coordinates), which slightly complicates the implementation of indirect method, but, importantly, does not affect their efficiency. \\

\noindent 2) \textit{Indirect method based on numerical homotopy procedures:} Our second aim consists of providing a numerical algorithm based on indirect methods. The main advantage of indirect methods is their extremely good numerical accuracy. Indeed, they inherit of the very quick convergence properties of the Newton method. Nevertheless, it is known that their main drawback is related to their initialization. We address this issue by adopting \textit{homotopy methods} (see, e.g., \cite{allgower2003introduction}).

The basic idea of homotopy methods is to solve a difficult problem step by step starting from a simpler problem (that we call \textit{problem of order zero}) by parameter deformation. Combined with the shooting problem derived from the PMP, homotopies consist of deforming the problem into a simpler one (i.e., on which a shooting method can be easily initialized) and then of solving a series of shooting problems step by step to come back to the original problem. One of the main issues is then being able to design an appropriate problem of order zero, which should ``resemble'' some extent of the initial problem but at the same time should be ``easy to solve''.

Homotopy procedures have proved to be reliable and robust for problems like orbit transfer, atmospheric reentry or planar tilting maneuvers (see, e.g., \cite{EmmanuelH,Petit,zhu2016minimum,zhu2016planar}). Here, we propose a numerical homotopy scheme to solve the shooting problem coming from the optimal guidance framework, ensuring high numerical accuracy of optimal trajectories.

To practically show the efficiency of this homotopy algorithm, we give numerical solutions of the \textit{endo-atmospheric missile interception} problem (presented, for example, in \cite{cottrell1971optimal}). We design an appropriate problem of order zero which is a good candidate to initialize the first homotopic iterations. Then, we solve the original problem by a \textit{linear continuation method} (i.e., the simplest homotopy scheme, see, e.g., \cite{allgower2003introduction}).

\subsection{Structure of the Paper} The paper is organized as follows. Section \ref{SectionProblem} contains details on the model under consideration and the optimal guidance problem. Section \ref{secTheoretical} is devoted to the PMP formulation of our problem, its intrinsic geometric behavior analysis and the computations of the optimal controls as functions of the state and the costate (which represents a crucial step to correctly define numerical indirect methods. Singular controls are analyzed too). In Sections \ref{sectHomotopy} and \ref{sectAppl} we provide the numerical homotopy scheme, giving global numerical solutions for the endo-atmospheric missile interception problem. Finally, Section \ref{conclSect} contains conclusions and perspectives.

\section{Optimal Guidance Problem} \label{SectionProblem}

\subsection{Model Dynamics for Guidance Systems} \label{SectionForces}

We focus on a class of launch vehicles modeled as a three-dimensional axial symmetric cylinder, where $\bm{u}$ denotes its principal body axis, steered by a control system (for example, based on steering fins or a reaction control system). We denote by $Q$ the point of the vehicle where this system is placed. Let $O$ be the center of the Earth, $\bm{K}$ be the northsouth axis of the planet and consider an orthonormal inertial frame $(\bm{I},\bm{J},\bm{K})$ centered at $O$. For the applications presented, the effect of the rotation of the Earth can be neglected. Denoting by $G$ the center of mass of the vehicle which is assumed to lie on $\bm{u}$, the motion is described by the variables $(\bm{r}(t),\bm{v}(t),\bm{u}(t))$, where $\bm{r}(t) = x(t) \bm{I} + y(t) \bm{J} + z(t) \bm{K}$ is the trajectory of $G$ while the vector $\bm{v}(t) = \dot{x}(t) \bm{I} + \dot{y}(t) \bm{J} + \dot{z}(t) \bm{K}$ denotes its velocity.

We denote by $m$ the mass of the vehicle, whose evolution is given as function of the mass flow rate, denoted by $q$. The air density is denoted by $\rho(\bm{r})$ (a standard exponential law of type $\rho_0 \exp(-(\|\bm{r}\|-r_T) / h_r)$ is considered, where $\rho_0 > 0$, $r_T$ is the radius of the Earth and $h_r$ is a reference altitude) while $S$ denotes a constant reference surface for aerodynamical forces. The forces that act on the vehicle are (see, e.g., \cite{pucci2015nonlinear,carlucci2018ballistics}):

\begin{itemize}
\item gravity $\bm{g} = -g(\bm{r}) m \frac{\bm{r}}{\| \bm{r} \|}$;
\item drag $\bm{D} = -\frac{S}{2} \rho(\bm{r}) C_D \| \bm{v} \| \bm{v}$, where \begingroup \small $C_D = C_{D_0} + C_{D_1} \left( \frac{\| \bm{u} \times \bm{v} \|}{\| \bm{v} \|} \right)^2$ \endgroup is a quadratic approximation of the drag coefficient ($C_{D_0}$, $C_{D_1}$ are constant;
\item lift $\bm{L} = \frac{S}{2} \rho(\bm{r}) C_{L_{\alpha}} \big( \bm{v} \times (\bm{u} \times \bm{v}) \big)$, where $C_{L_{\alpha}}$ is constant;
\item thrust $\bm{T} = f_T(t) \bm{u}$, where $f_T$ is a given nonnegative function which is proportional to the mass flow $q$.
\end{itemize}

Structural optimization ensures that torques do not affect the dynamics of the momentum. As a standard result (see, e.g. \cite{pucci2015nonlinear,carlucci2018ballistics}), the following dynamics is obtained

\begingroup
\small
\begin{eqnarray} \label{firstDyn}
\begin{cases}
\displaystyle \dot{\bm{r}}(t) = \bm{v}(t) , \ \dot{\bm{v}}(t) = \bm{f}(t,\bm{r}(t),\bm{v}(t),\bm{u}(t)) := \frac{\bm{T}(t,\bm{u}(t))}{m(t)} + \medskip \\
\displaystyle \hspace{10pt} \frac{\bm{g}(\bm{r}(t))}{m(t)} + \frac{\bm{D}(\bm{r}(t),\bm{v}(t),\bm{u}(t))}{m(t)} + \frac{\bm{L}(\bm{r}(t),\bm{v}(t),\bm{u}(t))}{m(t)} \ .
\end{cases}
\end{eqnarray}
\endgroup

\subsection{General Optimal Guidance Problem}

System \eqref{firstDyn} must be closed with some stability constraints. In particular, the velocity must be always positively oriented w.r.t. the principal body axis and, for controllability reasons, the velocity $\bm{v}$ must lies inside a cone whose symmetry axis is the body axis $\bm{u}$, and that has amplitude $0 < \alpha_{\max} \le \pi/6$, where $\alpha_{\max}$ is the \textit{maximal angle of attack}. From this and \eqref{firstDyn}, the full dynamics of our system becomes
\begingroup
\small
\begin{eqnarray} \label{guidanceDyn}
\begin{cases}
\dot{\bm{r}}(t) = \bm{v}(t) \quad , \quad \dot{\bm{v}}(t) = \bm{f}(t,\bm{r}(t),\bm{v}(t),\bm{u}(t)) \medskip \\
\displaystyle (\bm{r}(t),\bm{v}(t)) \in N \quad , \quad \bm{u}(t) \in S^2 \medskip \\
\bm{r}(0) = \bm{r}_0 \ , \ \bm{v}(0) = \bm{v}_0 \quad , \quad  (\bm{r}(T),\bm{v}(T)) \in M \subseteq N \medskip \\
\displaystyle c_1(\bm{v}(t),\bm{u}(t)) := -\bm{v}(t) \cdot \bm{u}(t) \le 0 \medskip \\
\displaystyle c_2(\bm{v}(t),\bm{u}(t)) := \bigg( \frac{\| \bm{u}(t) \times \bm{v}(t) \|}{\| \bm{v}(t) \| \sin \alpha_{\max}} \bigg)^2 - 1 \le 0
\end{cases}
\end{eqnarray}
\endgroup
where $N$ is an open subset of $\mathbb{R}^6 \setminus \{ 0 \}$ consisting of all possible scenarios (see Remark \ref{remN} in Section \ref{localPMP}), $S^2 = \{ \bm{u} \in \mathbb{R}^3 : \| \bm{u} \|^2 = 1 \}$ is the unit sphere in $\mathbb{R}^3$, $(\bm{r}_0,\bm{v}_0) \in N$ are given initial values, $T$ is the final time and $M$ is a subset of $N$ representing given final conditions. The control variable on which we act is represented by the principal body axis $\bm{u}$.

In this general context, a mission depends on which specific task the launch vehicle has to accomplish, which in turn depends on the cost that has to be minimized and on the set $M$ of final conditions. Then, given any function $g : \mathbb{R} \times \mathbb{R}^3 \times \mathbb{R}^3 \rightarrow \mathbb{R}$ of class $C^1$, we define the General Optimal Guidance Problem (\textbf{GOGP}) to be the optimal control problem that consists of minimizing the generic cost
\begin{equation*}
C(T,\bm{r}(\cdot),\bm{v}(\cdot),\bm{u}(\cdot)) = g(T,\bm{r}(T),\bm{v}(T))
\end{equation*}
under the dynamical control system (\ref{guidanceDyn}). The final time $T$ may be free or not. The generality of this cost allows one to consider various launch vehicle missions: for instance, in the case of endo-atmospheric landing problem one wants to minimize the error between the final position and some desired target point, or, in the case of missile interception one may want to maximize the final velocity.

In what follows, to apply indirect methods it will be needed to compute optimal controls using the PMP (see also Section \ref{controlSection}). This may become difficult to accomplish unless one considers further (merely technical) assumptions on $g$ and $M$. More specifically, we assume the following:
\begin{assumption} \label{assM}
The set $M$ is a submanifold of $N$. Moreover, at least one between the following two conditions is satisfied:
\begin{enumerate}
\item[A)] The final time $T$ is free and $\displaystyle \frac{\partial g}{\partial t}(T,\bm{r},\bm{v}) \neq 0$;

\item[B)] It holds $M = \Big\{ (\bm{r},\bm{v}) \in N : F(\bm{r},\bm{v}) = 0 \Big\}$, where $F$ is a smooth submersion. Moreover, for every local chart $(x_1,\dots,x_6)$ (local coordinates) of $(\bm{r},\bm{v}) \in M$ in $N$, there exists a variable $x_i$ such that $\frac{\partial g}{\partial x_i}(T,\bm{r},\bm{v}) \neq 0$.
\end{enumerate}
\end{assumption}

\section{Pontryagin Maximum Principle Analysis and Optimal Controls in Two Local Charts} \label{secTheoretical}

\subsection{PMP for Problems with Mixed Control-State Constraints} \label{maxSect}

The main objective of this paper consists of providing a numerical strategy to solve (\textbf{GOGP}) via indirect methods. They are based on necessary conditions for optimality that arise by applying the PMP to (\textbf{GOGP}) (see, e.g., \cite{betts1998survey}): in this section, we recall such necessary conditions for optimality.

The formulation of (\textbf{GOGP}) contains two mixed control-state constraints: $c_1$ and $c_2$. In the presence of such kind of constraints, the PMP can be efficiently employed only under further regularity assumptions on $c_1$ and $c_2$ (see, e.g., \cite{dmitruk2009development}). Indeed, it is required that the \textit{rank condition}
\begin{equation} \label{rankCond}
\textnormal{rank} \left( \begin{array}{cc}
\partial_{u_1} c_2 & u _1 \\
\partial_{u_2} c_2 & u_2 \\
\partial_{u_3} c_2 & u_3
\end{array} \right)(\bm{v},\bm{u}) = 2
\end{equation}
holds when $c_2(\bm{v},\bm{u}) = 0$ and $\| \bm{u} \|^2 = 1$ (see, e.g., \cite{dmitruk1993maximum,dmitruk2009development}). Straightforward computations show that \eqref{rankCond} is always satisfied, therefore, the PMP can be applied to (\textbf{GOGP}), leading to the following necessary conditions for optimality as follows.

Denote $\bm{p} = (\bm{p}_1,\bm{p}_2) \in \mathbb{R}^3 \times \mathbb{R}^3$, $\bm{\mu} = (\mu_0,\mu_1,\mu_2) \in \mathbb{R}^3$ and, as usual in the framework of the PMP, define
\begingroup
\begin{multline} \label{ham}
H(t,\bm{r},\bm{v},\bm{p},\bm{\mu},\bm{u}) := H^0(t,\bm{r},\bm{v},\bm{p},\bm{u}) + \mu_0 ( \| \bm{u} \|^2 - 1 )  \\
+ \mu_1 c_1(\bm{v},\bm{u}) + \mu_2 c_2(\bm{v},\bm{u}) := \Big( \bm{p}_1 \cdot \bm{v} + \bm{p}_2 \cdot \bm{f}(t,\bm{r},\bm{v},\bm{u}) \Big) \\
+ \mu_0 ( \| \bm{u} \|^2 - 1 ) + \mu_1 c_1(\bm{v},\bm{u}) + \mu_2 c_2(\bm{v},\bm{u})
\end{multline}
\endgroup
to be the \textit{Hamiltonian} of (\textbf{GOGP}) (see, e.g., \cite{dmitruk1993maximum,dmitruk2009development}). According to the PMP with mixed control-state constraints (see, e.g. \cite{pontryagin1987mathematical,hestenes1965calculus,dmitruk2009development}), if $(\bm{r}(\cdot),\bm{v}(\cdot),\bm{u}(\cdot))$ is optimal for (\textbf{GOGP}) with final time $T$, there exist a non-positive scalar $p^0$, an absolutely continuous curve $\bm{p} : [0,T] \rightarrow \mathbb{R}^6$ called \textit{adjoint vector}, and functions $\mu_0(\cdot)$, $\mu_1(\cdot)$, $\mu_2(\cdot) \in L^{\infty}([0,T],\mathbb{R})$, with $(\bm{p}(\cdot),p^0) \neq 0$, such that the so-called \textit{extremal} $(\bm{r}(\cdot),\bm{v}(\cdot),\bm{p}(\cdot),p^0,\mu_0(\cdot),\mu_1(\cdot),\mu_2(\cdot),\bm{u}(\cdot))$ satisfies almost everywhere in the time-interval $[0,T]$:

\begin{itemize}
\item \textbf{Adjoint Equations}
\begingroup
\begin{eqnarray} \label{adjointSystem}
\begin{cases}
\displaystyle \left(\begin{array}{c} \dot{\bm{r}}(t) \\ \dot{\bm{v}}(t) \end{array}\right) = \frac{\partial H}{\partial \bm{p}}(t,\bm{r}(t),\bm{v}(t),\bm{p}(t),\bm{\mu}(t),\bm{u}(t)) \medskip \\
\displaystyle \dot{\bm{p}}(t) = -\frac{\partial H}{\partial (\bm{r},\bm{v})}(t,\bm{r}(t),\bm{v}(t),\bm{p}(t),\bm{\mu}(t),\bm{u}(t))
\end{cases}
\end{eqnarray}
\endgroup
\item \textbf{Maximality Conditions}
\begingroup
\begin{equation} \label{maxCond}
\displaystyle H^0(t,\bm{r}(t),\bm{v}(t),\bm{p}(t),\bm{u}(t)) \ge H^0(t,\bm{r}(t),\bm{v}(t),\bm{p}(t),\bm{u})
\end{equation}
for every vector $\bm{u} \in S^2$ that satisfies $c_1(\bm{v}(t),\bm{u}) \le 0$ and $c_2(\bm{v}(t),\bm{u}) \le 0$. Moreover, it holds
\begin{equation} \label{maxCondDer}
\displaystyle \frac{\partial H}{\partial \bm{u}}(t,\bm{r}(t),\bm{v}(t),\bm{p}(t),\bm{\mu}(t),\bm{u}(t)) = 0
\end{equation}
\endgroup
\item \textbf{Complementarity Slackness Conditions}
\begingroup
\begin{equation} \label{slack}
\begin{cases}
\mu_1(t) c_1(\bm{v}(t),\bm{u}(t)) = 0 \\
\mu_2(t) c_2(\bm{v}(t),\bm{u}(t)) = 0
\end{cases} \ , \ \mu_1(t) \le 0 \ , \ \mu_2(t) \le 0
\end{equation}
\endgroup
\item \textbf{Transversality Conditions}
\begingroup
\begin{equation} \label{transvCond2}
\displaystyle \bm{p}(T) - p^0 \frac{\partial g}{\partial (\bm{r},\bm{v})}(T,\bm{r}(T),\bm{v}(T)) \perp T_{(\bm{r}(T),\bm{v}(T))} M
\end{equation}
\endgroup
where $T_{(\bm{r}(T),\bm{v}(T))} M$ is the tangent space of $M$ at $(\bm{r}(T),\bm{v}(T)) \in M$. If the final time $T$ is free, then
\begingroup
\small
\begin{equation} \label{transvCond1}
\displaystyle \max_{\bm{u}} H^0(T,\bm{r}(T),\bm{v}(T),\bm{p}(T),\bm{u}) = - p^0 \frac{\partial g}{\partial t}(T,\bm{r}(T),\bm{v}(T))
\end{equation}
\endgroup
where the maximum is taken over vectors $\bm{u} \in S^2$ that satisfy $c_1(\bm{v}(T),\bm{u}) \le 0$ and $c_2(\bm{v}(T),\bm{u}) \le 0$.
\end{itemize}
The extremal is said \textit{normal} if $p^0 \neq 0$ and, in this case, we set $p^0 = -1$. Otherwise, the extremal is said \textit{abnormal}.

As pointed out in the introduction, obtaining rigorous and useful information on the multipliers $\mu_1(\cdot)$, $\mu_2(\cdot)$ may be difficult, which consequently makes challenging applying indirect methods, as stated by the previous conditions, to (\textbf{GOGP}).

\subsection{Local Model with Respect to Two Local Charts} \label{localPMP}

A change of coordinates can be used to transform the mixed control-state constraints $c_1$ and $c_2$ into pure control constraints, allowing to use standard indirect methods. This is commonly used in aerospace (see, e.g., \cite{bonnard2003optimal}), though without the global (geometric) insight that we propose in this paper. However, this transformation acts only locally, preventing one from representing the whole configuration manifold $N$. For sake of clarity, we first recall this standard transformation, and then, we show how to fix the problem of Euler singularities by introducing further coordinates, in which, $c_1$ and $c_2$ remain pure control constraints. In turn, this allows us to locally apply standard indirect methods to solve (\textbf{GOGP}), by employing a simplified version of \eqref{adjointSystem}-\eqref{transvCond1}. \\

\subsubsection{Reduction to Pure Control Constraints via Local Coordinates} \label{firstChart}

We denote by $(r,L,\ell)$ the spherical coordinates of the center of mass $G$ of the vehicle w.r.t. $(\bm{I},\bm{J},\bm{K})$, where $r$ is the distance between $O$ and $G$ (Section \ref{SectionForces}), $L$ the latitude and $\ell$ the longitude. We denote $(\bm{e}_L,\bm{e}_{\ell},\bm{e}_r)$ the \textit{North-East-Down} (NED) frame, a moving frame centered at $G$, where $-\bm{e}_r$ is the local vertical direction, $(\bm{e}_L,\bm{e}_{\ell})$ is the local horizontal plane and $\bm{e}_L$ is pointing to the North. By definition, we have
\begingroup
\begin{eqnarray*}
\begin{cases}
\bm{e}_L = -\sin(L) \cos({\ell}) \bm{I} - \sin(L) \sin({\ell}) \bm{J} + \cos(L) \bm{K} \\
\bm{e}_{\ell} = -\sin({\ell}) \bm{I} + \cos({\ell}) \bm{J} \\
\bm{e}_r = -\cos(L) \cos({\ell}) \bm{I} - \cos(L) \sin({\ell}) \bm{J} - \sin(L) \bm{K}
\end{cases}
\end{eqnarray*}
\endgroup
for which $\bm{r} = -r \bm{e}_r$ and we have
\begingroup
\begin{equation} \label{derNED}
\begin{split}
&\dot{\bm{e}}_L = - \dot{\ell} \sin(L) \bm{e}_{\ell} + \dot{L} \bm{e}_r \ , \ \dot{\bm{e}}_{\ell} =  \dot{\ell} \sin(L) \bm{e}_L + \dot{\ell} \cos(L) \bm{e}_r \\
&\dot{\bm{e}}_r = -\dot{L} \bm{e}_L  - \dot{\ell} \cos(L) \bm{e}_{\ell} \ .
\end{split}
\end{equation}
\endgroup
Then, the transformation from the frame $(\bm{I},\bm{J},\bm{K})$ to the frame $(\bm{e}_L,\bm{e}_{\ell},\bm{e}_r)$ is the following rotation (i.e., a mapping in $SO(3) = \{ R \in GL_3(\mathbb{R}) : R^{\top} R = I \ , \ \textnormal{det}(R) = 1 \})$
\begingroup
\small
\begin{equation*}
R(L,{\ell}) := \left( \begin{array}{ccc}
-\sin(L) \cos({\ell}) & -\sin(L) \sin({\ell}) & \cos(L) \\
-\sin({\ell}) & \cos({\ell}) & 0 \\
-\cos(L) \cos({\ell}) & -\cos(L) \sin({\ell}) & -\sin(L) \end{array} \right) .
\end{equation*}
\endgroup
\begin{wrapfigure}{r}{0.2\textwidth}
\centering
\includegraphics[width=0.18\textwidth]{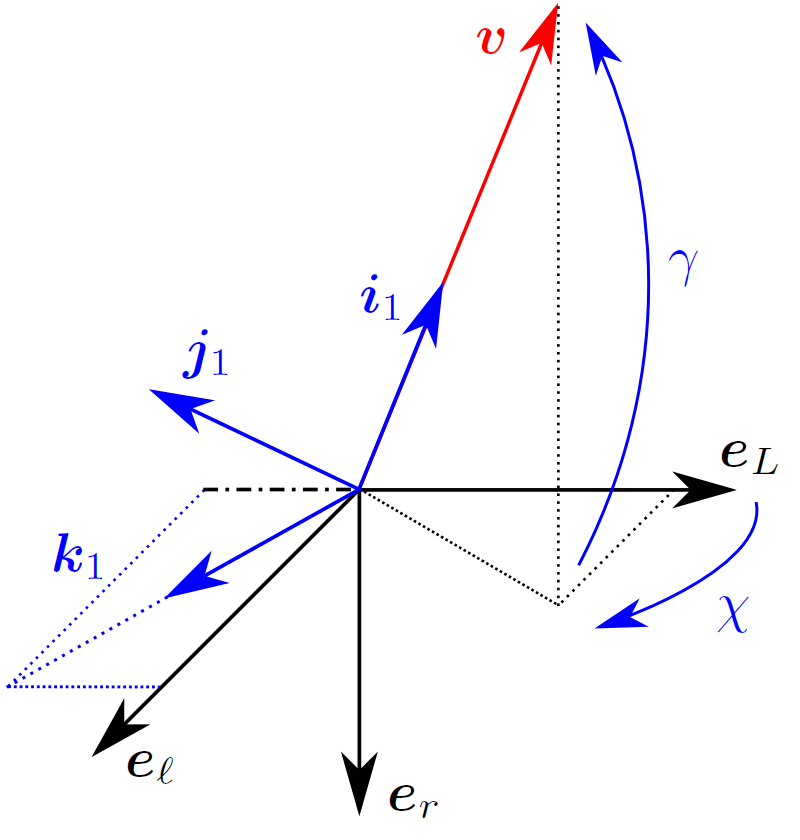}
\vspace{-10pt}
\caption{Frame $(\bm{i}_1,\bm{j}_1,\bm{k}_1)$.}
\end{wrapfigure}

To obtain $c_1$ and $c_2$ as pure control constraints, further coordinates for the velocity of the vehicle must be introduced. Using the classical formulation in the azimuth/path angle coordinates (see, e.g., \cite{bonnard2003optimal}), the \textit{first velocity frame} $(\bm{i}_1,\bm{j}_1,\bm{k}_1)$ is

\begingroup
\small
\begin{eqnarray} \label{frame1}
\begin{cases}
\displaystyle \bm{i}_1 := \frac{\bm{v}}{v} = \cos(\gamma) \cos(\chi) \bm{e}_L + \cos(\gamma) \sin(\chi) \bm{e}_{\ell} - \sin(\gamma) \bm{e}_r \medskip \\
\bm{j}_1 := -\sin(\gamma) \cos(\chi) \bm{e}_L - \sin(\gamma) \sin(\chi) \bm{e}_{\ell} - \cos(\gamma) \bm{e}_r \\
\bm{k}_1 := -\sin(\chi) \bm{e}_L + \cos(\chi) \bm{e}_{\ell}
\end{cases}
\end{eqnarray}
\endgroup
where we denote $v = \| \bm{v} \|$. Therefore, the rotation from the frame $(\bm{e}_L,\bm{e}_{\ell},\bm{e}_r)$ to the frame $(\bm{i}_1,\bm{j}_1,\bm{k}_1)$ is
\begingroup
\small
\begin{equation*}
R_a(\gamma,\chi) = \left( \begin{array}{ccc}
\cos(\gamma) \cos(\chi) & \cos(\gamma) \sin(\chi) & -\sin(\gamma) \\
-\sin(\gamma) \cos(\chi) & -\sin(\gamma) \sin(\chi) & -\cos(\gamma) \\
-\sin(\chi) & \cos(\chi) & 0 \end{array} \right) .
\end{equation*}
\endgroup

It is important to note that $(r,L,{\ell},v,\gamma,\chi)$ represent local coordinates for the dynamics of (\textbf{GOGP}). In the context of differential geometry, this means that there exists a local chart of $\mathbb{R}^6 \setminus \{ 0 \}$ whose coordinates are exactly $(r,L,{\ell},v,\gamma,\chi)$. Indeed, denote $U = \Big[ (0,\infty) \times \left( -\frac{\pi}{2} , \frac{\pi}{2} \right) \times ( -\pi , \pi ) \Big]^2$ and define the mapping $\varphi^{-1}_a : U \longrightarrow \mathbb{R}^6 \setminus \{ 0 \}$ such that
\begingroup
\begin{multline} \label{map1}
\varphi^{-1}_a(r,L,{\ell},v,\gamma,\chi) = \bigg( r \cos(L) \cos({\ell}) , r \cos(L) \sin({\ell}) , \\
r \sin(L) , R^{\top}(L,{\ell}) R^{\top}_a(\gamma,\chi) \Bigg( \begin{array}{c} v \\ 0 \\ 0 \end{array} \Bigg) \bigg) .
\end{multline}
\endgroup
This mapping is an injective immersion and its inverse is a local chart of $\mathbb{R}^6 \setminus \{ 0 \}$ (in the sense of differential geometry) when restricted to $U_a := \varphi^{-1}_a(U)$, which is an open subset of $\mathbb{R}^6 \setminus \{ 0 \}$. Exploiting (\ref{derNED}) and the definition of $(\bm{i}_1,\bm{j}_1,\bm{k}_1)$, in the coordinates provided by (\ref{map1}), the derivative of $\bm{v}$ is
\begingroup
\begin{multline} \label{velocity1}
\dot{\bm{v}} = \dot{v} \bm{i}_1 + \left( v \dot{\gamma} - \frac{v^2}{r} \cos(\gamma) \right) \bm{j}_1 + \\
\left( v \cos(\gamma) \dot{\chi} - \frac{v^2}{r} \cos^2(\gamma) \sin(\chi) \tan(L) \right) \bm{k}_1 \ .
\end{multline}
\endgroup

Finally, we introduce new control variables (which are functions of the original control $\bm{u}$), under which, $c_1$ and $c_2$ can be reformulated as pure control constraints. For this, define the new control $\bm{w} = R_a(\gamma,\chi) R(L,{\ell}) \bm{u}$. Then, the constraint functions become (by using the fact that $v > 0$ by definition)
\begingroup
\begin{equation} \label{localConstraint1}
c_1(\bm{w}) = -w_1 \ , \ c_2(\bm{w}) = \frac{w^2_2 + w^2_3}{\sin^2(\alpha_{\max})} - 1 \ , \ \bm{w} \in S^2
\end{equation}
\endgroup
which are pure control constraints. Denote the normalized drag and lift coefficients respectively by $d = \frac{1}{2 m} \rho S C_{D_0}$, $c_m = \frac{1}{2 m} \rho S C_{L_{\alpha}}$ and the efficiency factor by $\eta > 0$ (see, e.g., \cite{pucci2015nonlinear,carlucci2018ballistics,pepy2014indirect}). By introducing $\omega(t) = \frac{f_T(t)}{m(t) v(t)} + v(t) c_m(t) > 0$, with the help of (\ref{velocity1}), the local evaluation of the dynamics in (\ref{guidanceDyn}) using the local chart $\varphi_a$ immediately gives
\begingroup
\small
\begin{eqnarray} \label{dynFirst}
\begin{cases}
\dot{r} = v\sin(\gamma) \ , \ \dot{L} = \displaystyle \frac{v}{r} \cos(\gamma) \cos(\chi) \ , \ \dot{\ell} = \displaystyle \frac{v}{r} \frac{\cos(\gamma) \sin(\chi)}{\cos(L)} \medskip \\
\dot{v} = \displaystyle \frac{f_T}{m} w_1 -\left(d + \eta c_m (w^2_2 + w^2_3) \right) v^2 - g \sin(\gamma) \medskip \\
\dot{\gamma} = \displaystyle \omega w_2 + \left(\frac{v}{r} - \frac{g}{v}\right) \cos(\gamma) \medskip \\
\dot{\chi} = \displaystyle \frac{\omega}{\cos(\gamma)} w_3 + \frac{v}{r} \cos(\gamma) \sin(\chi) \tan(L) \ .
\end{cases}
\end{eqnarray}
\endgroup
It is crucial to note that $\gamma = \pm\pi/2$ are singularities for \eqref{dynFirst}.

The previous computations allow us to reformulate (\textbf{GOGP}) by introducing a new optimal control problem, named (\textbf{GOGP})$_a$, which si locally equivalent to (\textbf{GOGP}) and has only pure control constraints: this represents one of the two sought optimal control problems on which we run classical indirect methods. It consists of minimizing the cost
\begin{equation*}
C_a(T,r,L,{\ell},v,\gamma,\chi,\bm{w}) = g(T,\varphi^{-1}_a(r,L,{\ell},v,\gamma,\chi)(T))
\end{equation*}
subject to the dynamics (\ref{dynFirst}) and the control constraints (\ref{localConstraint1}). \\

\subsubsection{Additional Coordinates to Manage Eulerian Singularities} \label{secondChart}

Even if formulation (\textbf{GOGP})$_a$ is widely used in the aerospace community, it prevents one from completely describing the original problem (\textbf{GOGP}) because of its local nature. Indeed, in several situations, demanding performance criteria (costs $C$) and onerous missions (final conditions $M$) force optimal trajectories to pass through points that do not lie within the domain of the local chart $\varphi_a$ (i.e., $U_a$), and then, by exploiting merely (\textbf{GOGP})$_a$ either the optimality could be lost or, in the worst case, the numerical computations may fail.

\begin{wrapfigure}{r}{0.24\textwidth}
\vspace{-10pt}
\centering
\includegraphics[width=0.25\textwidth]{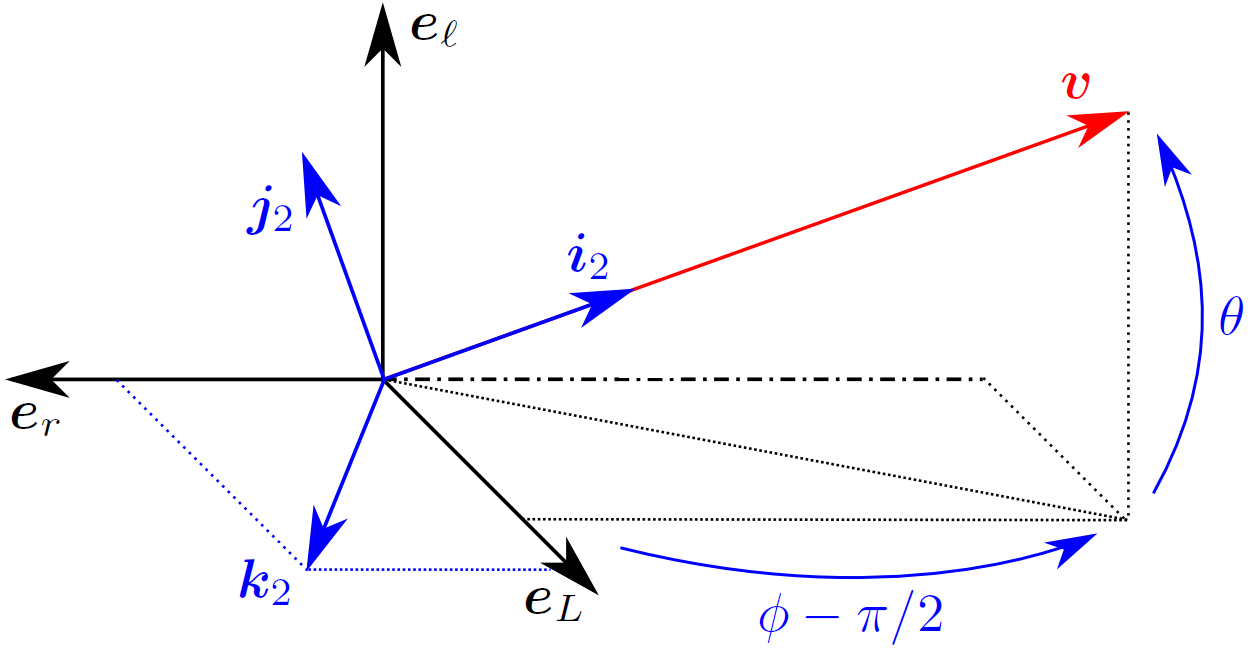}
\vspace{-20pt}
\caption{Frame $(\bm{i}_2,\bm{j}_2,\bm{k}_2)$.}
\vspace{-15pt}
\end{wrapfigure}
Here, the novelty consists of introducing another set of coordinates that covers the singularities (with respect to the path angle $\gamma$) of chart $(U_a,\varphi_a)$ in which the constraints $c_1$ and $c_2$ are pure control constraints, as provided by expressions (\ref{localConstraint1}).

For this, by mimicking the previous case, we introduce a new, \textit{second velocity frame} $(\bm{i}_2,\bm{j}_2,\bm{k}_2)$, defined as

\begingroup
\small
\begin{eqnarray} \label{frame2}
\begin{cases}
\displaystyle \bm{i}_2 = \frac{\bm{v}}{v} = \cos(\theta) \sin(\phi) \bm{e}_L + \sin(\theta) \bm{e}_{\ell} + \cos(\theta) \cos(\phi) \bm{e}_r \medskip \\
\bm{j}_2 = -\sin(\theta) \sin(\phi) \bm{e}_L + \cos(\theta) \bm{e}_{\ell} - \sin(\theta) \cos(\phi) \bm{e}_r \\
\bm{k}_2 = -\cos(\phi) \bm{e}_L + \sin(\phi) \bm{e}_r
\end{cases}
\end{eqnarray}
\endgroup
and the transformation (rotation) from the frame $(\bm{e}_L,\bm{e}_{\ell},\bm{e}_r)$ to the frame $(\bm{i}_2,\bm{j}_2,\bm{k}_2)$ is given by
\begingroup
\small
\begin{equation*}
R_b(\theta,\phi) = \left( \begin{array}{ccc}
\cos(\theta) \sin(\phi) & \sin(\theta) & \cos(\theta) \cos(\phi) \\
-\sin(\theta) \sin(\phi) & \cos(\theta) & -\sin(\theta) \cos(\phi) \\
-\cos(\phi) & 0 & \sin(\phi) \end{array} \right) .
\end{equation*}
\endgroup

The new local chart $(U_b,\varphi_b)$ is given by its domain $U_b = \varphi^{-1}_b(U)$ (see Section \ref{firstChart} for the definition of $U$) and
\begingroup
\begin{multline*}
\varphi^{-1}_b(r,L,{\ell},v,\theta,\phi) = \bigg( r \cos(L) \cos({\ell}) , r \cos(L) \sin({\ell}) , \\
r \sin(L) , R^{\top}(L,{\ell}) R^{\top}_b(\theta,\phi) \Bigg( \begin{array}{c} v \\ 0 \\ 0 \end{array} \Bigg) \bigg) .
\end{multline*}
\endgroup
This new local chart covers the singularities with respect to the path angle $\gamma$ of the local chart $(U_a,\varphi_a)$. In these new coordinates, the derivative of the velocity is
\begingroup
\begin{multline} \label{velocity2}
\dot{\bm{v}} = \dot{v} \bm{i}_2 + \bigg( v \dot{\theta} - \frac{v^2}{r} \sin(\theta) \big( \cos(\phi) + \sin(\phi) \tan(L) \big) \bigg) \bm{j}_2  \\
+ \bigg( \frac{v^2}{r} \cos^2(\theta) \Big( \sin(\phi) + \tan^2(\theta) \big( \sin(\phi) - \tan(L) \cos(\phi) \big) \Big) \\
+ v \dot{\phi} \cos(\theta) \bigg) \bm{k}_2 \ .
\end{multline}
\endgroup

As in the previous case, we now introduce new control variables (which are complementary to the local control $\bm{w}$), by defining $\bm{z} = R_b(\theta,\phi) R(L,{\ell}) \bm{u}$. Simple computations show that the constraints $c_1$ and $c_2$ are given in this local chart by
\begingroup
\begin{equation} \label{localConstraint2}
c_1(\bm{z}) = -z_1 \ , \ c_2(\bm{z}) = \frac{z^2_2 + z^2_3}{\sin^2(\alpha_{\max})} - 1 \ , \ \bm{z} \in S^2 \ .
\end{equation}
\endgroup
Using the same notations as in the previous case, with the help of expression (\ref{velocity2}), the local evaluation of the dynamics in (\ref{guidanceDyn}) by using the local chart $\varphi_b$ immediately gives
\begingroup
\footnotesize
\begin{eqnarray} \label{dynSecond}
\begin{cases}
\dot{r} = -v\cos(\theta)\cos(\phi) \ , \ \dot{L} = \displaystyle \frac{v}{r} \cos(\theta) \sin(\phi) \ , \ \dot{\ell} = \displaystyle \frac{v}{r} \frac{\sin(\theta)}{\cos(L)} \medskip \\
\dot{v} = \displaystyle \frac{f_T}{m} z_1 -\left(d + \eta c_m (z^2_2 + z^2_3) \right) v^2 + g \cos(\theta) \cos(\phi) \medskip \\
\dot{\theta} = \displaystyle \omega z_2 + \frac{v}{r} \sin(\theta) \Big( \cos(\phi) + \sin(\phi) \tan(L) \Big) - \frac{g}{v} \sin(\theta) \cos(\phi) \medskip \\
\dot{\phi} = \displaystyle -\frac{\omega}{\cos(\theta)} z_3 + \frac{v}{r} \cos(\theta) \bigg( \sin(\phi) + \tan^2(\theta) \Big( \sin(\phi) \medskip \\
\displaystyle \hspace{5ex} - \tan(L) \cos(\phi) \Big) \bigg) - \frac{g}{v} \frac{\sin(\phi)}{\cos(\theta)} \ .
\end{cases}
\end{eqnarray}
\endgroup

We define a second optimal control problem, named (\textbf{GOGP})$_b$, which is locally equivalent to (\textbf{GOGP}) and has only pure control constraints: this represents the second sought optimal control problem on which we run classical indirect methods. It consists of minimizing the cost
\begin{equation*}
C_b(T,r,L,{\ell},v,\theta,\phi,\bm{z}) = g(T,\varphi^{-1}_b(r,L,{\ell},v,\theta,\phi)(T))
\end{equation*}
subject to the dynamics (\ref{dynSecond}) and the control constraints (\ref{localConstraint2}).

\begin{remark} \label{remN}
The mappings $\varphi^{-1}_a : U \rightarrow \mathbb{R}^6 \setminus \{ 0 \}$, $\varphi^{-1}_b : U \rightarrow \mathbb{R}^6 \setminus \{ 0 \}$ are not defined respectively for the values $\chi = \pi$, $\phi = \pi$: these singularities can be covered by extending $\varphi^{-1}_a$ and $\varphi^{-1}_b$ also on $\Big[ (0,\infty) \times \left( -\frac{\pi}{2} , \frac{\pi}{2} \right) \times ( 0 , 2 \pi ) \Big]^2$. Nevertheless, the framework of this paper concerns launch vehicles able to cover bounded distances (in the region of one hundred kilometers). From these remarks, without loss of generality, we define the configuration manifold of (\textbf{GOGP}) to be $N := U_a \cup U_b$.
\end{remark}

\subsection{Equivalence between Global and Local Formulations} \label{secChange}

From the previous sections, it is clear that, within the open set $U_a \subseteq \mathbb{R}^6 \setminus \{ 0 \}$, (\textbf{GOGP}) is equivalent to (\textbf{GOGP})$_a$ while, within the open set $U_b \subseteq \mathbb{R}^6 \setminus \{ 0 \}$, (\textbf{GOGP}) is equivalent to (\textbf{GOGP})$_b$. However, it is not clear whether the PMP formulation related to (\textbf{GOGP}), which is a problem with mixed control-state constraints, is equivalent respectively to the dual formulation of (\textbf{GOGP})$_a$ locally within $U_a$, and with the dual formulation of (\textbf{GOGP})$_b$ locally within $U_b$, which are problems with pure control constraints. More precisely, we have a priori three different tuples of multipliers, namely:
\begin{itemize}
\item $(\bm{p}(\cdot),p^0,\mu_1(\cdot),\mu_2(\cdot))$ related to (\textbf{GOGP});
\item $(p_a(\cdot),p^0_a)$ related to (\textbf{GOGP})$_a$;
\item $(p_b(\cdot),p^0_b)$ related to (\textbf{GOGP})$_b$.
\end{itemize}
Nothing ensures that $p_a(\cdot)$, $p_b(\cdot)$ are related to $\bm{p}(\cdot)$ within $U_a$, $U_b$, respectively. In such situation, recasting the analysis of necessary conditions for optimality and related indirect methods from (\textbf{GOGP}) to (\textbf{GOGP})$_a$, (\textbf{GOGP})$_b$ may cause inconsistencies because, a priori, the adjoint vectors $p_a(\cdot)$, $p_b(\cdot)$, coming from the local formulations, evolve independently.

We fix this gap by showing that $p_a(\cdot)$, $p_b(\cdot)$ can be consistently related to $\bm{p}(\cdot)$, which will justify the study and the development of indirect methods for (\textbf{GOGP})$_a$, (\textbf{GOGP})$_b$ to solve (\textbf{GOGP}). In particular, we prove that it is always possible to choose the previous multipliers so that the local projections of $(\bm{p}(\cdot),p^0)$ onto charts $(U_a,\varphi_a)$ and $(U_b,\varphi_b)$ are equivalent respectively to $(p_a(\cdot),p^0_a)$ and to $(p_b(\cdot),p^0_b)$.

\medskip

\begin{theorem} \label{mainTheo}
\textit{Consider the manifold $N = U_a \cup U_b \subseteq \mathbb{R}^6 \setminus \{ 0 \}$ of all possible scenarios for (\textbf{GOGP}). Suppose that $(\bm{r}(\cdot),\bm{v}(\cdot),\bm{u}(\cdot))$ is an optimal solution for (\textbf{GOGP}) in $[0,T]$. There exist multipliers $(\bm{p}(\cdot),p^0,\mu_1(\cdot),\mu_2(\cdot))$ satisfying the PMP formulation (\ref{adjointSystem})-(\ref{transvCond1}) and multipliers $(p_a(\cdot),p^0_a)$, $(p_b(\cdot),p^0_b)$ related to the classical PMP formulations with pure control constraints respectively of problem (\textbf{GOGP})$_a$ and of problem (\textbf{GOGP})$_b$, such that $p^0_a = p^0_b = p^0$ and}
\begin{equation} \label{adjPull}
\bm{p}(t) = \begin{cases} (\varphi_a)^*_{\varphi_a(\bm{r}(t),\bm{v}(t))} \cdot p_a(t) \quad , \quad (\bm{r}(t),\bm{v}(t)) \in U_a \\
(\varphi_b)^*_{\varphi_b(\bm{r}(t),\bm{v}(t))} \cdot p_b(t) \quad , \quad (\bm{r}(t),\bm{v}(t)) \in U_b \end{cases}
\end{equation}
\textit{where $(\cdot)^*$ is the pullback (see \cite{agrachev2013control} for such definitions).}
\end{theorem}

\medskip

The proof of Theorem \ref{mainTheo} is done in Appendix \ref{proofApp}. The main idea is the following. By the PMP for problems with mixed control-state constraints, there exists a global adjoint vector $\bm{p}(\cdot)$ for (\textbf{GOGP}) which we restrict it to the domain of one of the two local charts built previously, for instance, $(U_a,\varphi_a)$. Then, via the local maximality condition (\ref{maxCondDer}) and the transformation between $\bm{u}$ and $\bm{w}$ (see Sections \ref{firstChart}, \ref{secondChart}), one shows that the covector $(\varphi^{-1}_a)^* \cdot \bm{p}(\cdot)$ satisfies the PMP formulation with pure control constraints related to (\textbf{GOGP})$_a$.

Let us clarify how one can exploit this result to solve (\textbf{GOGP}) by indirect methods. Assume to have an optimal solution $(\bm{r}(\cdot),$ $\bm{v}(\cdot),$ $\bm{u}(\cdot))$ for (\textbf{GOGP}) in $[0,T]$. Without loss of generality, we can assume that $(\bm{r},\bm{v})(0) \in U_a$. If a guess for the optimal value of $p_a(0)$ (or equivalently of $\bm{p}(0)$, see \eqref{adjPull}) is known, we can solve (\textbf{GOGP}) by running an indirect method on (\textbf{GOGP})$_a$ starting from $p_a(0)$. Suppose that, at a given time $\tau_1 \in (0,T)$, the optimal trajectory is such that $(\bm{r},\bm{v})(\tau_1) \in U_b \setminus U_a$, i.e., our solution crosses a \textit{singular region} of the first local chart. Then, Theorem \ref{mainTheo} allows us to stop the numerical computations at a time $\tau_2 < \tau_1$ such that $(\bm{r},\bm{v})(\tau_2) \in U_a \cap U_b$ and then run an indirect method on (\textbf{GOGP})$_b$ starting from $p_b(\tau_2) = (\varphi_a \circ \varphi^{-1}_b)^*_{\varphi_a(\bm{r}(\tau_2),\bm{v}(\tau_2))} p_a(\tau_2)$ (see \eqref{adjPull}), therefore avoiding the singularity related to $\varphi_a$ when reaching the point $(\bm{r},\bm{v})(\tau_1) \in U_b \setminus U_a$ (see Figure \ref{fig:charts} below). This procedure can be iterated every time a jump from $U_a$ to $U_b$ (as well as a jump from $U_b$ to $U_a$) occurs in the optimal trajectory. The adjoint vector related to (\textbf{GOGP}) is recovered thanks to (\ref{adjPull}). This methodology allows one to describe global optimal solutions for any feasible mission related to problem (\textbf{GOGP}).

It is worth noting that, even if modifying indirect methods by implementing the change of coordinates \eqref{adjPull} introduces further computations, the numerical transformation between local adjoint vectors takes a negligible part of the total computational time an indirect method needs to converge (as simulations show, see Section \ref{secNumerical}), justifying our approach.

\begin{figure}[htpb]
\centering
\includegraphics[width=0.25\textwidth]{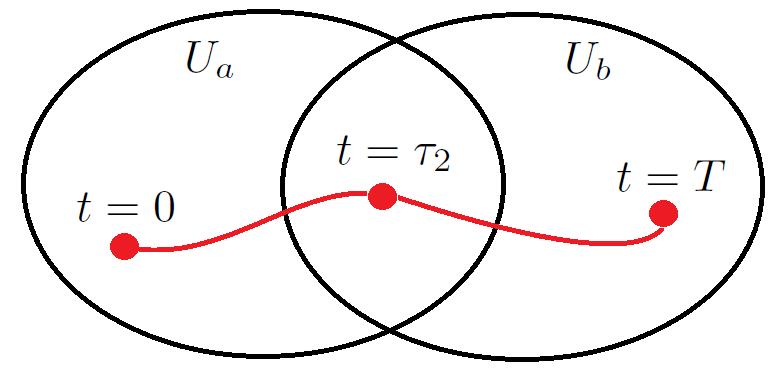}
\caption{Optimal trajectory crossing the domains of the two local charts.} \label{fig:charts}
\end{figure}

\subsection{Optimal Control as Functions of the Adjoint Vectors} \label{controlSection}

In the previous section, we showed (by Theorem \ref{mainTheo}) that implementing indirect methods on (\textbf{GOGP}) is equivalent to developing indirect methods for the local PMP formulations related to (\textbf{GOGP})$_a$ and to (\textbf{GOGP})$_b$. In this section, we provide optimal controls as functions of the states and the adjoint vectors, by adopting the previous formalism based on (\textbf{GOGP})$_a$, (\textbf{GOGP})$_b$. This provides formulations needed to run indirect methods on problems (\textbf{GOGP})$_a$, (\textbf{GOGP})$_b$. \\

Let $(\bm{r}(\cdot),\bm{v}(\cdot),\bm{u}(\cdot))$ be an optimal solution for (\textbf{GOGP}) in $[0,T]$, and $\bm{p}(\cdot)$, $\bm{p}_a(\cdot) = (p^a_r,p^a_L,p^a_{\ell},p^a_v,p_{\gamma},p_{\chi})(\cdot)$ and $\bm{p}_b(\cdot) = (p^b_r,p^b_L,p^b_{\ell},p^b_v,p_{\theta},p_{\phi})(\cdot)$ be the related adjoint vectors respectively for (\textbf{GOGP}), (\textbf{GOGP})$_a$ and for (\textbf{GOGP})$_b$ as in Theorem \ref{mainTheo} (see also Section \ref{localPMP} for the definition of the local problems). As pointed out previously, thanks to Theorem \ref{mainTheo}, the computation of the optimal control $\bm{u}$ can be achieved by focusing on the optimal values of the local controls $\bm{w}$, $\bm{z}$, which are the projections of $\bm{u}$ onto $U_a$, $U_b$, respectively (see also Section \ref{localPMP}). Hereafter, when clear from the context, we skip the dependence on $t$ to keep better readability. For sake of clarity, we denote $C_a := p^a_v \frac{f_T}{m}$, $C_b := p^b_v \frac{f_T}{m}$, $D_a := p^a_v \eta c_m v^2$ and $D_b := p^b_v \eta c_m v^2$ (see Section \ref{localPMP} for notations).

Expressions for optimal controls $\bm{w}$, $\bm{z}$ as functions of the local states and the local adjoint vectors can be achieved by studying the local versions of the Maximality Condition \eqref{maxCond}. From the PMP for pure control constraints applied to (\textbf{GOGP})$_a$, (\textbf{GOGP})$_b$ (resulting as a special case of conditions \eqref{adjointSystem}-\eqref{transvCond1}, see also \cite{pontryagin1987mathematical}), locally almost everywhere where they are defined, respectively related to (\textbf{GOGP})$_a$ and to (\textbf{GOGP})$_b$, these maximality conditions are given by

\begingroup
\small
\begin{multline} \label{firstHam}
\displaystyle \bm{w}(t) = \argmax \Bigg\{ C_a w_1 - D_a (w^2_2 + w^2_3) + p_{\gamma} \omega w_2 + p_{\chi} \frac{\omega}{\cos(\gamma)} w_3 \\
w^2_1 + w^2_2 + w^2_3 = 1 \ , \ w_1 \ge 0 \ , \ w^2_2 + w^2_3 \le \sin^2(\alpha_{\max}) \Bigg\}
\end{multline}
\endgroup

\begingroup
\small
\begin{multline} \label{secondHam}
\displaystyle \bm{z}(t) = \argmax \Bigg\{ C_b z_1 - D_b (z^2_2 + z^2_3) + p_{\theta} \omega z_2 - p_{\phi} \frac{\omega}{\cos(\theta)} z_3 \\
z^2_1 + z^2_2 + z^2_3 = 1 \ , \ z_1 \ge 0 \ , \ z^2_2 + z^2_3 \le \sin^2(\alpha_{\max}) \Bigg\} \ .
\end{multline}
\endgroup

Solving these maximization conditions may lead to either \textit{regular} or \textit{nonregular controls}, depending on the value of $(p_{\gamma}(\cdot),p_{\chi}(\cdot))$, $(p_{\theta}(\cdot),p_{\phi}(\cdot))$ on non-zero measure subsets of $[0,T]$. Indeed, by definition, regular controls are the regular points of \textit{the end-point mapping} while nonregular controls are its critical points (see, e.g., \cite{trelat2012optimal}). Then, for (\textbf{GOGP}), regular controls consist of controls whose extremal satisfies either $p_{\gamma}|_J(\cdot) \neq 0$ or $p_{\chi}|_J(\cdot) \neq 0$, within a non-zero measure subset $J \subseteq [0,T]$, if the system travels along the first local chart $(U_a,\varphi_a)$ within $J$. On the other hand, regular controls satisfy $p_{\theta}|_J(\cdot) \neq 0$ or $p_{\phi}|_J(\cdot) \neq 0$ if the system covers the second local chart $(U_b,\varphi_b)$ within $J$. Conversely, nonregular controls consist of controls for which there exists a non-zero measure subset $J \subseteq [0,T]$ such that $p_{\gamma}|_J(\cdot) = p_{\chi}|_J(\cdot) = 0$ in the first chart, and $p_{\theta}|_J(\cdot) = p_{\phi}|_J(\cdot) = 0$ in the second chart.

We analyze separately regular and nonregular controls. \\

\subsubsection{Regular Controls} Suppose that, locally within a non-zero measure subset $J \subseteq [0,T]$, either $p_{\gamma}|_J(\cdot) \neq 0$ or $p_{\chi}|_J(\cdot) \neq 0$ if the system travels along the first chart $(U_a,\varphi_a)$ within $J$. Otherwise, $p_{\theta}|_J(\cdot) \neq 0$ or $p_{\phi}|_J(\cdot) \neq 0$. In this case, regular controls appear. Explicit expressions for these are easily derived from (\ref{firstHam}), (\ref{secondHam}) by using the Karush-Kuhn-Tucker conditions (see, e.g., \cite{NoceWrig06}), if we assume the following:

\medskip

\begin{assumption} \label{assumpSimply}
\textit{For points $(\varepsilon,x) \in \mathbb{R}_+ \times \mathbb{R}$ such that $(1 + \varepsilon) x^2 \le \sin^2(\alpha_{\max})$, where $0 < \alpha_{\max} \le \pi / 6$ is constant, the following holds: $\sqrt{1 - (1 + \varepsilon) x^2} \cong \Big( 1 - (1 + \varepsilon) x^2 / 2 \Big)$.}
\end{assumption}

\medskip

This assumption is not limiting. Indeed, for most of the applications that are based on the dynamical model developed for (\textbf{GOGP}), the maximal angle of attack $\alpha_{\max}$ is actually lower than $\pi / 6$ (because of controllability issues). Moreover, this assumption has already implicitly been used to recover the explicit expressions of the drag and the lift listed in Section \ref{SectionForces}. Under Assumption \ref{assumpSimply}, we provide the computations for the explicit expressions of regular controls in Appendix \ref{regularApp}. Note that regular controls are well defined in each of the two charts $(U_a,\varphi_a)$, $(U_b,\varphi_b)$, but their local expressions reach singular values as soon as the optimal trajectory gets close to the boundary of $U_a$, $U_b$, respectively. \\

\subsubsection{Nonregular Controls} \label{secSingular}

In some cases, locally within a non-zero measure subset $J \subseteq [0,T]$, it may happen that $p_{\gamma}|_J(\cdot) = p_{\chi}|_J(\cdot) = 0$ in the first local chart, or $p_{\theta}|_J(\cdot) = p_{\phi}|_J(\cdot) = 0$ in the second local chart. The control is then nonregular and the evaluation of optimal controls is harder to achieve than in the regular case. Here, (\ref{firstHam}) and (\ref{secondHam}) are

\begingroup
\begin{multline} \label{singularHam}
\displaystyle \bm{w}(t) = \argmax \Big\{ C_a w_1 - D_a (w^2_2 + w^2_3) \, \mid \,
w^2_1 + w^2_2 + w^2_3 = 1 , \\
w_1 \ge 0 , w^2_2 + w^2_3 \le \sin^2(\alpha_{\max}) \Big\}
\end{multline}
\endgroup

\begingroup
\begin{multline} \label{singularHam2}
\displaystyle \bm{z}(t) = \argmax \Big\{ C_b z_1 - D_b (z^2_2 + z^2_3) \ \mid \ z^2_1 + z^2_2 + z^2_3 = 1 , \\
z_1 \ge 0 \ , \ z^2_2 + z^2_3 \le \sin^2(\alpha_{\max}) \Big\} \ .
\end{multline}
\endgroup

The Karush-Kuhn-Tucker conditions are no more helpful because, depending on the value of $C_a$ or $C_b$, many uncountable values of $(w_2,w_3)$ or $(z_2,z_3)$ are optimal. Instead, a geometric study is required. It is in the case of nonregular controls that Assumption \ref{assM} becomes particularly useful to manage hard computations, as well as the following one:

\medskip

\begin{assumption} \label{assVel}
\textit{Let $J \subseteq [0,T]$ be a non-zero measure subset. Along $J$, any optimal trajectory associated with a nonregular control in $J$ satisfies (see Section \ref{localPMP} for notations)}
\begingroup
$$
\| \bm{v} \|^2 > \frac{3}{2} g(\bm{r}) h_r \bigg( \sqrt{1 + \frac{4}{9} \frac{1}{g(\bm{r}) h_r} \left( \frac{f_T}{m d} \right) } - 1 \bigg) \ .
$$
\endgroup
\end{assumption}

\medskip

It is important to note that, for our applications, the magnitude of the velocities of the vehicles is in general large enough when $f_T > 0$, so that Assumption \ref{assVel} is always satisfied, as numerical simulations confirm. In particular, this assumption is required only for nonregular arcs, i.e., in the case of regular optimal controls no boundaries on the velocities are imposed. Under Assumption \ref{assVel}, we provide the expressions of nonregular optimal controls in Appendix \ref{singularApp}, which, together with the results above, lead straightforwardly to the following:

\medskip

\begin{proposition} \label{regSing}
\textit{Under Assumption \ref{assumpSimply}, regular optimal controls for (\textbf{GOGP}) are well-defined and have univocal explicit expressions. On the other hand, under Assumption \ref{assM} and Assumption \ref{assVel}, any nonregular optimal control for (\textbf{GOGP}) is well-defined and has a univocal explicit expression.}
\end{proposition}

\medskip

It is worth noting that nonregular controls are not very common: the absence of nonregular controls has been widely studied in the context of optimal control and it has been shown that, under appropriate assumptions, regular controls appear \textit{almost always} (see, e.g., \cite{bonnard1997generic,chitour2008singular}). Running Monte-Carlo simulations on (\textbf{GOGP}) for many different realistic scenarios in the context of missile interception problems, we have never found nonregular controls (see Section \ref{secBatch}). However, for sake of completeness, we have provided full descriptions of both regular and nonregular controls so that indirect methods for (\textbf{GOGP}) are always correctly defined (see Proposition \ref{regSing}).

\section{Numerical Indirect Method for General Optimal Guidance Problems} \label{sectHomotopy}

In the previous sections, we showed that (\textbf{GOGP}) can be locally converted into two optimal control problems with pure control constraints, and that the PMP formulations of such problems locally match with the original global one. This allows one to run indirect methods on these local problems to provide solutions for (\textbf{GOGP}) with transmission conditions. In this section, we describe a general numerical scheme to solve (\textbf{GOGP}) which combines shooting methods with homotopy procedures (see, e.g., \cite{allgower2003introduction}). The main idea consists of introducing a family of problems parametrized by some quantity $\bm{\lambda}$, so that: 1) the problem related to $\bm{\lambda} = 0$, named \textit{problem of order zero}, is simple to solve by shooting methods; 2) we solve (\textbf{GOGP}) by combining shooting methods with an iterative procedure that makes $\bm{\lambda}$ vary with continuity, starting from the solution obtained for $\bm{\lambda} = 0$. We begin the section by discussing the formulation of the problem of order zero.

\subsection{Designing the Problem of Order Zero} \label{subSect0}

The problem of order zero, from which the iterative shooting path starts, should be, on one hand, handy to solve via basic shooting methods and, on the other hand, as close as possible to (\textbf{GOGP}) to efficiently recover a solution for the original problem by some homotopy. It is worth noting that two difficulties prevent shooting methods to easily converge: 1) the presence of the gravity, the thrust and Earth's curvature that considerably complexify the dynamics; 2) the more demanding the mission is (which is represented by the cost $g$ and the target $M$), the less intuition one has on the structure of optimal solutions. These facts led us to derive the following heuristic problem of order zero for (\textbf{GOGP}): denoted by (\textbf{GOGP})$_0$, it consists of minimizing the simplified cost
\begin{equation*}
C_0(T,\bm{r}(\cdot),\bm{v}(\cdot),\bm{u}(\cdot)) = g_0(T,\bm{r}(T),\bm{v}(T))
\end{equation*}
subject to the simplified dynamics
\begingroup
\small
\begin{eqnarray*}
\begin{cases}
\dot{\bm{r}}(t) = \bm{v}(t) \quad , \quad \dot{\bm{v}}(t) = \bm{f}_0(t,\bm{r}(t),\bm{v}(t),\bm{u}(t)) \medskip \\
(\bm{r}(t),\bm{v}(t)) \in N \quad, \quad \bm{u}(t) \in S^2 \medskip \\
\displaystyle \bm{r}(0) = \bm{r}_0 \ , \ \bm{v}(0) = \bm{v}_0 \quad , \quad (\bm{r}(T),\bm{v}(T)) \in M_0 \subseteq N \medskip \\
\displaystyle c_1(\bm{v}(t),\bm{u}(t)) \le 0 \quad , \quad c_2(\bm{v}(t),\bm{u}(t)) \le 0 \ .
\end{cases}
\end{eqnarray*}
\endgroup
Here, $g_0$, $\bm{f}_0$ and $M_0$ represent simplified versions for the original cost $g$, the original dynamics $f$, and the original target $M$, respectively. Consistently with our remarks above, we define the simplified dynamics $\bm{f}_0$ such that the contributions of the gravity, the thrust and Earth's curvature are removed from the original dynamics, that is
\begingroup
\small
\begin{equation} \label{simplyDyn}
\bm{f}_0(t,\bm{r},\bm{v},\bm{u}) := \bm{f}(t,\bm{r},\bm{v},\bm{u}) - \bigg( \frac{\bm{T}(t,\bm{u})}{m} + \frac{\bm{g}(\bm{r})}{m} - \bm{\omega}_{\textnormal{NED}}(\bm{r},\bm{v}) \times \bm{v} \bigg)
\end{equation}
\endgroup
where $\bm{\omega}_{\textnormal{NED}}(\bm{r},\bm{v})$ represents the angular velocity of the NED frame $(\bm{e}_L,\bm{e}_l,\bm{e}_r)$ w.r.t. the inertial frame $(\bm{I},\bm{J},\bm{K})$ (the sign minus in front of it is consistent with our convention, see Section \ref{firstChart}). It is important to evaluate (\ref{simplyDyn}) only by using charts $(U_a,\varphi_a)$, $(U_b,\varphi_b)$, otherwise its explicit expression could appear more complex than the original dynamics, especially because of the presence of the term $\bm{\omega}_{\textnormal{NED}}\times \bm{v}$.

From what we pointed out above, the simplified cost $g_0$ and target $M_0$ should be designed such that it is easy to make a shooting method converge for (\textbf{GOGP})$_0$. This may require that, for instance, we choose $g_0$, $M_0$ so that optimal trajectories for (\textbf{GOGP})$_0$ do not meet Euler singularities, i.e., they lie entirely within the chart domain $U_a$ (or $U_b$), or also, that optimal strategies are not of bang-bang type. In a very general context, it may be not evident to provide appropriate $g_0$, $M_0$ such that (\textbf{GOGP})$_0$ is easy to solve by shooting methods, and this strongly depends on the nature of the original mission. The engineer intuition is often crucial at this step. Also, exploiting techniques from geometric control or dynamical system theory applied to mission design (see, e.g., \cite{trelat2012optimal}) may lead to design relevant formulations for $g_0$, $M_0$. Hereafter, we show how to define (\textbf{GOGP})$_0$ for missile interception (see Section \ref{sectAppl}).

Solving (\textbf{GOGP})$_0$ by standard shooting methods leads to a solution $(\bm{r}_0(\cdot),\bm{v}_0(\cdot),\bm{u}_0(\cdot))$ for (\textbf{GOGP})$_0$ with adjoint variables $(\bm{p}_0(\cdot),p^0_0)$. Thanks to Theorem \ref{mainTheo}, from now on, we do not report the multipliers related to the mixed constraints.

\subsection{Homotopies Initialized by the Problem of Order Zero} \label{sectHomotopy2}

Now that the problem of order zero has been defined, we provide optimal strategies for (\textbf{GOGP}) by iteratively solving a sequence of shootings problems indexed by some parameter $\bm{\lambda}$, using the adjoint variables related to (\textbf{GOGP})$_0$.

We first introduce the family of problems, denoted by (\textbf{GOGP})$_{\bm{\lambda}}$, depending on the parameter $\bm{\lambda}$. For every $\bm{\lambda} = (\lambda_1,\lambda_2) \in [0,1]^2$, the optimal control problem (\textbf{GOGP})$_{\bm{\lambda}}$ consists of minimizing the parametrized cost
\begin{equation*}
C_{\bm{\lambda}}(T,\bm{r}(\cdot),\bm{v}(\cdot),\bm{u}(\cdot)) = g_{\bm{\lambda}}(T,\bm{r}(T),\bm{v}(T))
\end{equation*}
subject to the parametrized dynamics
\begingroup
\small
\begin{eqnarray*}
\begin{cases}
\dot{\bm{r}}(t) = \bm{v}(t) \quad , \quad \dot{\bm{v}}(t) = \bm{f}_{\bm{\lambda}}(t,\bm{r}(t),\bm{v}(t),\bm{u}(t)) \medskip \\
(\bm{r}(t),\bm{v}(t)) \in N \quad, \quad \bm{u}(t) \in S^2 \medskip \\
\displaystyle \bm{r}(0) = \bm{r}_0 \ , \ \bm{v}(0) = \bm{v}_0 \quad , \quad (\bm{r}(T),\bm{v}(T)) \in M_{\bm{\lambda}} \subseteq N \medskip \\
\displaystyle c_1(\bm{v}(t),\bm{u}(t)) \le 0 \quad , \quad c_2(\bm{v}(t),\bm{u}(t)) \le 0 \ .
\end{cases}
\end{eqnarray*}
\endgroup
Here, the parametrized cost $g_{\bm{\lambda}}$ and dynamics $\bm{f}_{\bm{\lambda}}$ are
\begingroup
\small
\begin{multline*}
g_{\bm{\lambda}}(T,\bm{r},\bm{v}) := g_0(T,\bm{r},\bm{v}) + \lambda_1 \Big( g(T,\bm{r},\bm{v}) - g_0(T,\bm{r},\bm{v}) \Big) \\
\bm{f}_{\bm{\lambda}}(t,\bm{r},\bm{v},\bm{u}) := \bm{f}_0(t,\bm{r},\bm{v},\bm{u}) + \lambda_1 \Big( \bm{f}(t,\bm{r},\bm{v},\bm{u}) - \bm{f}_0(t,\bm{r},\bm{v},\bm{u}) \Big)
\end{multline*}
\endgroup
involving only the first component of $\bm{\lambda}$, and the parametrized target $M_{\bm{\lambda}}$ is chosen such that it depends only on the second component of $\bm{\lambda}$ and satisfies $M_{\lambda_2 = 0} \equiv M_0$, $M_{\lambda_2 = 1} \equiv M$. Remark that the problem of order zero (\textbf{GOGP})$_0$ corresponds to $\bm{\lambda} = (0,0)$ while the original (target) problem (\textbf{GOGP}) corresponds to $\bm{\lambda} = (1,1)$. In this case, the homotopy procedure consists of solving a series of shooting problems by making $\bm{\lambda}$ continuously pass from $(0,0)$ to $(1,1)$.

We propose to solve (\textbf{GOGP}) by Algorithm \ref{ref:algoCont} below, continuation scheme which operates on the family of problems (\textbf{GOGP})$_{\bm{\lambda}}$ introduced above, starting from problem (\textbf{GOGP})$_0$.

\SetKwInOut{Input}{Input}
\SetKwInOut{Output}{Output}
\SetKwInOut{Data}{Data}

\begin{algorithm}
	\caption{Numerical Continuation for (\textbf{GOGP})} \label{ref:algoCont}
	\Input{Solution $(\bm{r}_0(\cdot),\bm{v}_0(\cdot),\bm{u}_0(\cdot))$ for (\textbf{GOGP})$_0$ with adjoint variables $(\bm{p}_0(\cdot),p^0_0)$.}
	\Output{Solution $(\bm{r}_{\bm{\lambda}}(\cdot),\bm{v}_{\bm{\lambda}}(\cdot),\bm{u}_{\bm{\lambda}}(\cdot))$ for (\textbf{GOGP})$_{\bm{\lambda}}$.}
	\Data{Maximal number of iterations $k_{\max}$.}
	
	\Begin
	{
		$\lambda_1 = 0$, $\lambda_2 = 0$, $\Delta = 1$, $k = 0$\\
		\For{$i = 1$, $i \leftarrow i + 1$, $i \le 2$}
		{
			\While{$\lambda_i < 1$ and $k \le k_{\max}$}{
				$\lambda^{\textnormal{temp}}_i = \lambda_i + \Delta$, \begingroup\small$\bm{\lambda}^{\textnormal{temp}} = \left\{\begin{array}{cc}
				(\lambda^{\textnormal{temp}}_1,\lambda_2) , & i = 1 \\
				(\lambda_1,\lambda^{\textnormal{temp}}_2) , & i = 2
				\end{array}\right.$\endgroup\\
				Solve (\textbf{GOGP})$_{\bm{\lambda}^{\textnormal{temp}}}$ by a shooting initialized with $(\bm{p}_{\bm{\lambda}}(\cdot),p^0_{\bm{\lambda}})$ related to (\textbf{GOGP})$_{\bm{\lambda}}$\\
				\If{succesful}{
					$\bm{\lambda} = \bm{\lambda}^{\textnormal{temp}}$
				}
				\Else{
					$\Delta = \Delta / 2$
				}
			}
			$\Delta = 1$, $k = 0$\\
		}
		\Return{$(\bm{r}_{\bm{\lambda}}(\cdot),\bm{v}_{\bm{\lambda}}(\cdot),\bm{u}_{\bm{\lambda}}(\cdot))$}
	}
\end{algorithm}

Algorithm \ref{ref:algoCont} operates a continuation via bisection method on the coordinates of parameter $\bm{\lambda}$ and is considered successful if it ends with $\bm{\lambda} = (1,1)$. If line \textbf{10} of Algorithm \ref{ref:algoCont} is called frequently, the convergence rate may become slow. To prevent such behavior, acceleration steps may be considered (see, e.g., \cite{allgower2003introduction}). Numerical simulations show that splitting the continuation on the hard terms of the dynamics and on the mission helps obtaining better performance. The theoretical convergence of Algorithm \ref{ref:algoCont} is established under appropriate assumptions. Indeed, it is known that homotopy methods may fail whenever bifurcation points, singularities or different connected components are found (see, e.g., \cite{allgower2003introduction,trelat2012optimal}). However, the absence of \textit{conjugate points} and of \textit{abnormal minimizers} (remark that these may be different from nonregular controls) are sufficient conditions for homotopies to converge (see \cite{trelat2012optimal}). Algorithm \ref{ref:algoCont} may be combined with numerical procedures computing conjugate points (see \cite{bonnard2007second}). For what concerns (\textbf{GOGP}), we have solved many different realistic scenarios via Monte-Carlo simulations for missile interception problems, showing an empirical efficiency for Algorithm \ref{ref:algoCont} in such context (in particular, see Section \ref{secIPOPT}).

\section{Launch Vehicle Application: Endo-Atmospheric Missile Interception} \label{sectAppl}

In this section, we apply Algorithm \ref{ref:algoCont} to solve (\textbf{GOGP}) in the context of \textit{endo-atmospheric interception} (see, e.g., \cite{cottrell1971optimal}). The problem consists of steering a missile towards a given target, optimizing some criterion. We are interested in the \textit{mid-course phase} which starts when the vehicle reaches a given threshold of the magnitude of the velocity. The target consists of a predicted interception point and, since this point may change over time, fast and accurate computations are needed.

The Optimal Interception Problem (\textbf{OIP}) consists of the specific (\textbf{GOGP}) for which the cost $g$ and the target $M$ are
\begingroup
\begin{equation} \label{interception}
\displaystyle g(T,\bm{r}(T),\bm{v}(T)) = C_1 T - \| \bm{v}(T) \|^2
\end{equation}
\endgroup
\begingroup
\begin{multline} \label{finalSubMissile}
\displaystyle M = \bigg\{ (\bm{r},\bm{v}) \in N \ \mid \ \bm{r} = \bm{r}_1 \ , \ \bm{v} \cdot \bm{e}_r = \| \bm{v} \| \cos(\psi_1) \ , \\
\bm{v} \cdot \bm{e}_L = \| \bm{v} \| \cos(\psi_2) \ , \ \bm{v} \cdot \bm{e}_l = \| \bm{v} \| \sin(\psi_2) \bigg\}
\end{multline}
\endgroup
where $C_1 \ge 0$ and the final time $T$ is either fixed or free. This cost is set up to maximize the chances to complete the mission with reasonable delays. Moreover, $M$ fixes the final position and orientation of the vehicle. Assumption \ref{assM} is satisfied.

For numerical simulations, a \textit{solid-fuel propelled missile} is employed, with the following numerical values:
\begin{itemize}
\item \begingroup \small$c_m(0) = 7.5 \cdot 10^{-4} \textnormal{m}^{-1}$, $d(0) = 5 \cdot 10^{-5} \textnormal{m}^{-1}$, $\eta = 0.442$, $h_r = 7500 \textnormal{m}$, $\alpha_{\max} = \pi/6$, $q_0 = 0.025 \textnormal{s}^{-1}$, $f^0_T = 37.5 \textnormal{m} \cdot \textnormal{s}^{-2}$\endgroup
\item \begingroup \small$\displaystyle \frac{q}{m(0)}(t) = \begin{cases}
q_0 \ , \ t \le 20 \\
0 \ , \ t > 20
\end{cases} , \frac{f_T}{m(0)}(t) = \begin{cases}
f^0_T \ , \ t \le 20 \\
0 \ , \ t > 20 \ .
\end{cases}$\endgroup
\end{itemize}
The shooting problems in Algorithm \ref{ref:algoCont} are solved using a C++ environment and \textit{hybrd.c} \cite{minPACK} while a fixed time-step explicit fourth-order Runge-Kutta method is used to integrate differential equations (whose number of integration steps varies between 250 and 350). Computations are done on a system Ubuntu 12.04 (32-bit), with 7.00 Gb of RAM.

\subsection{Choices for the Problem of Order Zero} \label{secOrderZero}

Without loss of generality, a problem of order zero (\textbf{OIP})$_0$ for (\textbf{OIP}) can be chosen such that its optimal trajectory lies in the domain of the local chart $(U_a,\varphi_a)$. The following problem of order zero is considered (see \cite{bonalli2017analytical})
\begingroup
\small
\begin{eqnarray*}
\begin{cases}
\qquad \min \quad -v^2(T) \quad , \quad (w_2,w_3) \in \mathbb{R}^2 \bigskip \\
\dot{r} = v \sin(\gamma) \ , \ \dot{L} = \displaystyle \frac{v}{r} \cos(\gamma) \cos(\chi) \ , \ \dot{l} = \displaystyle \frac{v}{r} \frac{\cos(\gamma) \sin(\chi)}{\cos(L)} \medskip \\
\dot{v} = -(d + \eta c_m (w^2_2 + w^2_3)) v^2 \ , \ \dot{\gamma} = v c_m w_2 \ , \ \dot{\chi} = \displaystyle \frac{v c_m}{\cos(\gamma)} w_3
\end{cases}
\end{eqnarray*}
\endgroup
where, consistently with the arguments in Section \ref{subSect0}, the contribution of the thrust, the gravity and of $\bm{\omega}_{\textnormal{NED}} \times \bm{v}$ are removed, no constraints on the controls are imposed and $C_1 = 0$. More specifically, the target set $M_0$ can be chosen such that (\textbf{OIP})$_0$ is feasible and the PMP applied to (\textbf{OIP})$_0$ allows one to recover an approximated explicit guidance law which successfully initializes shooting methods for (\textbf{OIP})$_0$. For the sake of conciseness, we do not report all details for such computations and the interested reader is referred to \cite{bonalli2017analytical}.

\begin{figure}[t!]
\hspace*{-45pt} \includegraphics[width=0.65\textwidth]{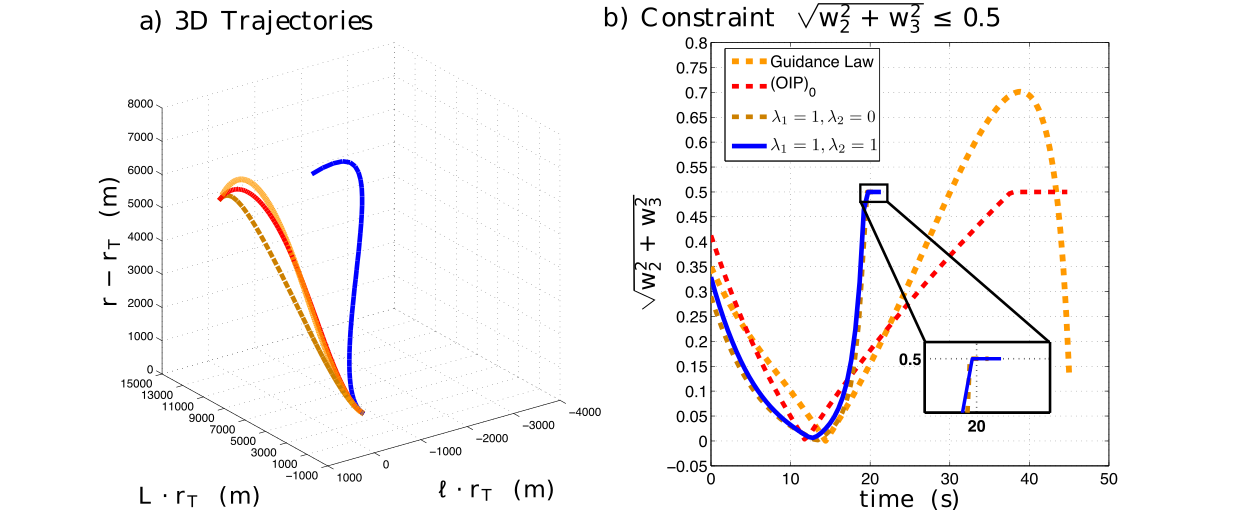}
\caption{Optimal trajectories $\bm{r}(\cdot)$ and constraints $c_2(\bm{r}(\cdot),\bm{v}(\cdot))$ (in the chart $(U_a,\varphi_a)$) for the first mission, i.e., cost $T - \| \bm{v}(T) \|^2$ and target $M^1$. The dashed-yellow curves represent the explicit guidance law used to initialize a shooting on (\textbf{OIP})$_0$ (see Section \ref{secOrderZero}), whose solution is dashed-red. The dashed-brown curves represent optimal quantities for $\bm{\lambda} = (1,0)$ and the solid-blue curves for $\bm{\lambda} = (1,1)$, i.e., for the original problem.}  \label{figMission1}
\end{figure}

\begin{figure}[t!]
\hspace*{-45pt} \includegraphics[width=0.65\textwidth]{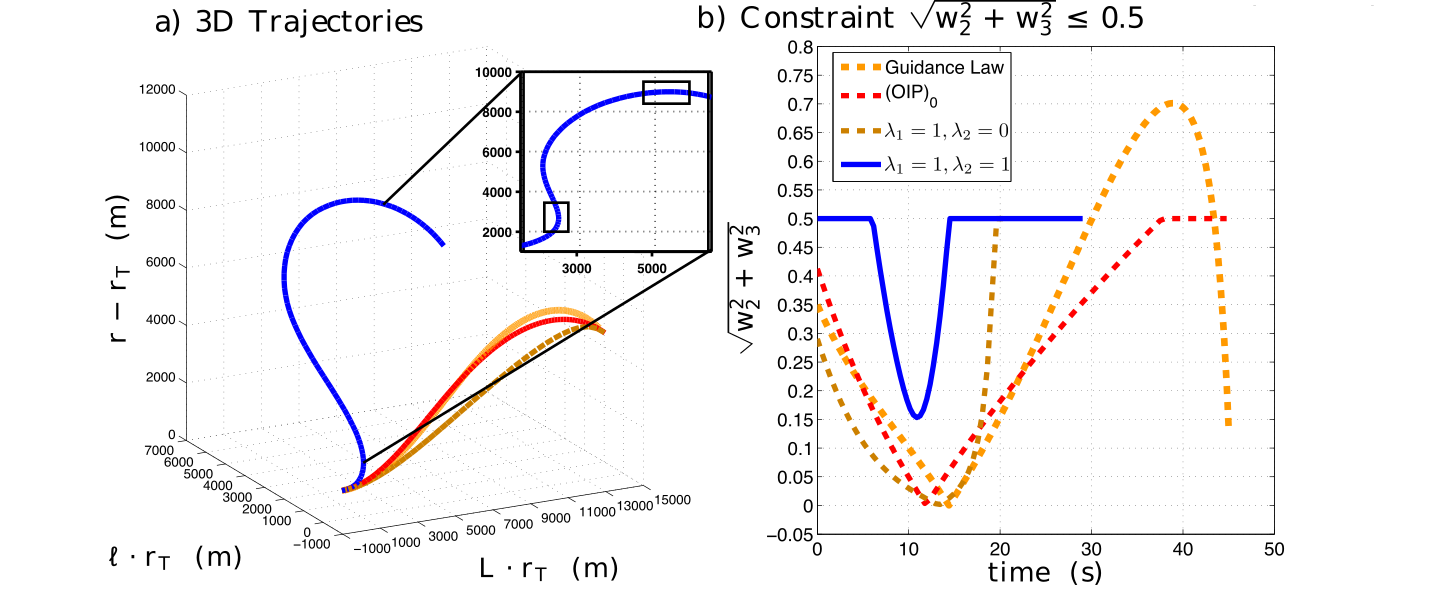}
\caption{Optimal trajectories $\bm{r}(\cdot)$ and constraints $c_2(\bm{r}(\cdot),\bm{v}(\cdot))$ (in the chart $(U_a,\varphi_a)$) for the second mission, i.e., cost $T - \| \bm{v}(T) \|^2$ and target $M^2$. The dashed-yellow curves represent the explicit guidance law used to initialize a shooting on (\textbf{OIP})$_0$ (see Section \ref{secOrderZero}), whose solution is dashed-red. The dashed-brown curves represent optimal quantities for $\bm{\lambda} = (1,0)$ and the solid-blue curves for $\bm{\lambda} = (1,1)$, i.e., for the original problem. In the big black box, we show the two-dimensional projection of the trajectory for (\textbf{OIP}) onto the plane $(L \cdot r_T,r-r_T)$. The small boxes show changes of local chart.}  \label{figMission2}
\end{figure}

\subsection{Numerical Simulations on Two Realistic Missions} \label{secNumerical}

In this section, we apply Algorithm \ref{ref:algoCont} to solve two realistic interception missions for (\textbf{OIP}). We consider free final time problems for which $C_1 = 1$ in \eqref{interception}. Details for the choice of the family of problems have been provided in Section \ref{sectHomotopy2}. In particular, the initial conditions and the target for the problem of order zero in standard units are (in the local chart $(U_a,\varphi_a)$)
$$
(r-r_T,L \cdot r_T,l \cdot r_T,v,\gamma,\chi)(0) = (1000,0,0,500,0,0)
$$
\begingroup
\small
$$
M_0 = \Big\{ (r-r_T,L \cdot r_T,l \cdot r_T) = (5000,14000,0) \ , \ (\gamma,\chi) = (-\pi/6,0) \Big\}
$$
\endgroup
where $r_T$ is Earth's radius. Since we fix \eqref{interception}, the two missions are unambiguously defined by \eqref{finalSubMissile}, i.e., respectively by
\begin{multline*}
M^1 = \Big\{ (r-r_T,L \cdot r_T,l \cdot r_T) = (5000,14000,-2000) \ , \\
(\gamma,\chi) = (-\pi/6,\pi/6) \Big\}
\end{multline*}
\begin{multline*}
M^2 = \Big\{ (r-r_T,L \cdot r_T,l \cdot r_T) = (7900,7500,2000) \ , \\
(\gamma,\chi) = (-\pi/4,-\pi/4) \Big\}
\end{multline*}
that are provided in standard units. The second mission is more challenging because abrupter maneuvers will be required. The parametrized target sets $M_{\lambda_2}$ for (\textbf{OIP})$_{\bm{\lambda}}$ are convex combinations in $\lambda_2$ of $M_0$ and $M^1$, $M^2$, respectively. \\

\subsubsection{First Mission}

When solving the first mission, Algorithm \ref{ref:algoCont} provides $(T,\| \bm{v}(T) \|) = (21.4,753.7)$ as optimal values (in standard units). Optimal trajectories and values for constraints are given in Figure \ref{figMission1}. The computations take around 1.6 s, for which 14 iterations on $\lambda_1$ and 11 iterations on $\lambda_2$ are required. This is due to the minimal time in the cost which makes the structure of the solutions more complex. In Figure \ref{figMission1} b), we see that, even if constraint $c_2$ is not satisfied by the explicit guidance law of Section \ref{secOrderZero}, the latter correctly initializes Algorithm \ref{ref:algoCont}, so that saturations on constraint $c_2$ are satisfied by the solution of the original problem. \\

\subsubsection{Second Mission}

For the second mission, considered to be more challenging, Algorithm \ref{ref:algoCont} provides $(T,\| \bm{v}(T) \|) = (29.03,475.2)$ as optimal values (in standard units). Optimal trajectories and values for constraints are given in Figure \ref{figMission2}. Here, the computations take around 2.3 s, where 14 iterations on $\lambda_1$ and 26 iterations on $\lambda_2$ are required. A higher number of iterations on $\lambda_2$ occurs because, when proposing to intercept a target quite close to the initial point, the vehicle is led to perform abrupt maneuvers to recover an optimal solution. Moreover, two changes of local chart as designed in Section \ref{secChange} are involved in this mission. Indeed, from Figure \ref{figMission2} a), we see that trajectories are close both to $\gamma = \pi/2$, singular value for $(U_a,\varphi_a)$, and to $\theta = \pi/2$, critical value for $(U_b,\varphi_b)$ (see the small black boxes in the two-dimensional projection onto $(L \cdot r_T,r-r_T)$ in Figure \ref{figMission2} a)). Even if the change of coordinates is not compulsory to solve this mission, it considerably increases the performances of the algorithm. Indeed, without it, simulations on this mission would take around 8.2 s with 19 iterations on $\lambda_1$ and 121 iterations on $\lambda_2$. Other tests show that some scenarios cannot be solved without the change of local chart. This shows a glimpse of the benefits in performance that one may achieve when adopting the change of local chart provided in this paper.

\subsection{Performance Test via Batch Simulations} \label{secBatch}

In this section, we test the efficiency of Algorithm \ref{ref:algoCont} on a batch of simulations. We consider free final time problems (\textbf{OIP}) without minimal final time, i.e., cost \eqref{interception} is chosen such that $C_1 = 0$. Monte-Carlo simulations are run on missions with the same initial conditions used in the previous simulations and values for the final target uniformly, randomly chosen in
$$
r - r_T \in [4000,8000] \ , \ L \cdot r_T \in [14000,18000] \ ,
$$
$$
l \cdot r_T \in [-4000,4000] \ , \ \gamma , \chi \in [-\pi/3,\pi/3]
$$
which gather realistic interception mission scenarios. We report in Table 1 the results obtained on 1000 missions.

Since the parametrized costs $C_{\lambda}$ do not change, the number of iterations on $\lambda_1$ is constant and equal to 7. The absence of minimal final time in the cost considerably improves computational times. Moreover, as mentioned at the end of Section \ref{secSingular}, every solution provided by Algorithm \ref{ref:algoCont} on this batch consists of regular controls, i.e., nonregular controls do not appear. From the fact that 99.7\% of the missions are solved, these tests empirically show that the structure of the problem of order zero that we described in Sections \ref{subSect0}, \ref{secOrderZero} for missile interception applications preserves enough information so that shooting methods combined with homotopies efficiently converge to a solution for the original problem. Remark that, without considering the change of local chart, an additional 9.5\% of missions would have failed.

\newcolumntype{C}{>{\centering\arraybackslash}p{1.85cm}}

\begingroup
\begin{table}
\centering
\footnotesize

\begin{tabular}{| C | C | C | C |}
\hline
Successful missions & Average time & Average nb. of iterations on $\lambda_2$ & Changes of local chart \\

\hline

99.7\% & 0.85 seconds & 6 & 9.5\% \\

\hline
\end{tabular}

\medskip

Table 1 : Average results obtained from solving (\textbf{OIP}) by Algorithm \ref{ref:algoCont} on 1000 missions. The averages are computed only on successful cases.
\end{table}
\endgroup

\subsection{Comparisons with State-of-the-Art Direct Methods} \label{secIPOPT}

In this section, we compare results and performance of Algorithm \ref{ref:algoCont} with state-of-the-art direct methods for optimization problem. We discretize (\textbf{OIP}) in time and solve the related nonlinear optimization problem in AMPL \cite{fourer1993ampl} combined with IpOpt \cite{wachter2006implementation} as nonlinear solver (see, e.g., \cite{gollmann2014theory,rodrigues2017optimal,ma2018direct}).

To have specific and fair results to compare, we run Algorithm \ref{ref:algoCont} and AMPL/IpOpt on three missions, by considering fixed final time problems (\textbf{OIP}) (and with $C_1 = 0$ in \eqref{interception}). The initial conditions are the same as in Sections \ref{secNumerical}, \ref{secBatch}, so that the missions are given by the following scenarios:
\begin{itemize}
\item (\textbf{SC})$_1$ : $T = 20$, $r - r_T =  7030$, $L \cdot r_T = 9450$, $l \cdot r_T = 1400$, $\gamma, \chi = -0.55$
\item (\textbf{SC})$_2$ : $T = 23$, $r - r_T =  7465$, $L \cdot r_T = 8475$, $l \cdot r_T = 1700$, $\gamma, \chi = -0.67$
\item (\textbf{SC})$_3$ : $T = 29$, $r - r_T =  7900$, $L \cdot r_T = 7500$, $l \cdot r_T = 2000$, $\gamma, \chi = -0.79$
\end{itemize}
in standard units. These missions are listed in ascending order by difficulty, so that the last mission requires a change of chart to be solved by Algorithm \ref{ref:algoCont}. The AMPL/IpOpt code is run considering the local chart $(U_a,\varphi_a)$, which is the formulation for (\textbf{OIP}) that is found in the literature. Moreover, to obtain comparable computational times, a second-order Runge-Kutta method is used to integrate ODEs both in Algorithm \ref{ref:algoCont} and in AMPL/IpOpt. We test various time-steps for AMPL/IpOpt, while setting 80 time-steps for Algorithm \ref{ref:algoCont}.

Results are given in Table 2. The finer the step-size of the time-discretization scheme for IpOpt is, the better optimal solutions are obtained, but this is at the price of additional computational time, which is particularly higher for the last, more challenging mission. For comparable computational times, we obtain slightly better solutions than AMPL/IpOpt.

\section{Conclusions and Perspectives} \label{conclSect}

In this paper we have developed a geometric analysis that we have used to design a numerical algorithm, based on indirect methods, to solve optimal control problems for endo-atmospheric launch vehicle systems. Considering the original problem with mixed control-state constraints in an intrinsic geometric framework, we have recast it into an optimal control problem with pure control constraints by restriction to two sets of local coordinates (local charts). We have solved the original problem by combining classical shooting methods and homotopies, bypassing singularities of Euler coordinates. We have provided numerical simulations for optimal interception missions, showing similar (sometimes better) performance than state-of-the-art methods in numerical optimal control.

Additional contributions may be considered to refine the dynamics. In particular, state and control delays may be important to take into account further dynamical strains and phenomena like the \textit{non-minimum phase}, a classical issue for launch vehicles applications (see, e.g., \cite{balas2012adaptive}). Motivated by the convergence result established in \cite{bonalli2017delay,bonalli2018delay}, we could add one further homotopic step on the delay. For computational times, even if many simulations on different missions for (\textbf{OIP}) show that Algorithm \ref{ref:algoCont} can run between 0.5 Hz and 1 Hz, we cannot ensure a real-time processing yet. This may be achieved by combining Algorithm \ref{ref:algoCont} with offline computations: we evaluate offline optimal strategies for several possible missions that will initialize online spatial continuations (i.e. on the parameter $\lambda_2$) to solve any new feasible mission.

\newcolumntype{D}{>{\centering\arraybackslash}p{0.75cm}}

\newcolumntype{E}{>{\centering\arraybackslash}p{1.5cm}}

\newcolumntype{F}{>{\centering\arraybackslash}p{2cm}}

\begingroup
\begin{table}
\centering
\footnotesize

\begin{tabular}{D | D | D | D | D | D | D |}

\cline{2-7}

& \multicolumn{2}{|F|}{(\textbf{SC})$_1$} & \multicolumn{2}{|F|}{(\textbf{SC})$_2$} & \multicolumn{2}{|F|}{(\textbf{SC})$_3$} \\

\cline{2-7}

& $\| v(T) \|$ & time & $\| v(T) \|$ & time & $\| v(T) \|$ & time \\

\hline

\multicolumn{1}{|E|}{Algorithm \ref{ref:algoCont}} & 763.1 & 0.66 & 609.8 & 0.92 & 480.0 & 2.3 \\

\hline

\multicolumn{1}{|E|}{IpOpt - 20 time steps} & 846.7 & 0.42 & 663.3 & 0.29 & 568.6 & 0.76 \\

\hline

\multicolumn{1}{|E|}{IpOpt - 40 time steps} & 780.3 & 1.1 & 635.7 & 1.6 & 531.2 & 1.9 \\

\hline

\multicolumn{1}{|E|}{IpOpt - 80 time steps} & 760.5 & 5.3 & 618.1 & 5.8 & 494.9 & 6.7 \\

\hline

\multicolumn{1}{|E|}{IpOpt - 120 time steps}  & 767.5 & 10 & 614.5 & 10 & 478.3 & 12 \\

\hline
\end{tabular}

\medskip

Table 2 : Optimal results obtained from solving (\textbf{OIP}) by Algorithm \ref{ref:algoCont} and AMPL/IpOpt. Quantities are reported in standard units.
\end{table}
\endgroup

\bibliographystyle{unsrt}
\bibliography{references}

\appendix

\subsection{Proof of the Consistency for Local Adjoint Vectors} \label{proofApp}

Here, we provide a proof of Theorem \ref{mainTheo}. By similarity between the local charts $(U_a,\varphi_a)$, $(U_b,\varphi_b)$, without loss of generality, we prove the assert considering the chart $(U_a,\varphi_a)$. \\

Denote $\bm{q} = (\bm{r},\bm{v})$ and, for the sake of clarity in the notation, let us denote an optimal solution for (\textbf{GOGP}) in $[0,T]$ by $(\bar{\bm{q}}(\cdot),\bar{\bm{u}}(\cdot))$. Select times $s_1, s_2 \in (0,T)$ such that $s_1 < s_2$ and $\bar{\bm{q}}([s_1,s_2]) \subseteq U_a$, and consider the notation $x = \varphi_a(\bm{q})$, $\bm{q} \in U_a$. Since $\bar{\bm{q}}([s_1,s_2]) \subseteq U_a$, in what follows, we merely need to consider (\textbf{GOGP}) when restricted to $U_a$. Therefore, when restricted to $[s_1,s_2]$, the curve $(\bar{\bm{q}}(\cdot),\bar{\bm{u}}(\cdot))$ is a solution for the following local version of (\textbf{GOGP}) in $U_a$ (which is correctly well-defined up to multiplications by appropriate smooth cut-off functions):
\begin{eqnarray*}
(\textbf{GOGP})_{\textnormal{loc}} \
\begin{cases}
\displaystyle \ \min \ g\big(s_2,\bm{q}(t),\bm{u}(t)\big) \medskip \\
\dot{\bm{q}}(t) = \bm{h}\big(t,\bm{q}(t),\bm{u}(t)\big) \ , \ \bm{q}(t) \in U_a \medskip \\
\bm{q}(s_1) = \bar{\bm{q}}(s_1) \ , \ \bm{q}(s_2) = \bar{\bm{q}}(s_2) \medskip \\
\bm{c}\big(\bm{q}(t),\bm{u}(t)\big) \le 0 \ , \ \textnormal{a.e.} \ [0,T]
\end{cases}
\end{eqnarray*}
where $s_1$, $s_2$ are fixed and, for sake of brevity, we denoted $\bm{h} = (\bm{v},\bm{f})$, $c_0(\bm{q},\bm{u}) := \| \bm{u} \|^2 - 1$ (even if $c_0$ does not depend on $\bm{q}$, this notation will be useful hereafter) and $\bm{c} = (c_0,c_1,c_2)$. On the other hand, with the definitions in Section \ref{localPMP} (see \eqref{localConstraint1}, \eqref{localConstraint2}), the local version of (\textbf{GOGP})$_{\textnormal{loc}}$ w.r.t. $(U_a,\varphi_a)$ writes as
\begin{eqnarray*}
(\textbf{GOGP})_a
\begin{cases}
\displaystyle \ \min \ g\big(s_2,\varphi_a^{-1}\big(x(t)\big),\Phi\big(x(t),\bm{w}(t)\big)\big) \medskip \\
\dot{x}(t) = d\varphi_a \cdot \bm{h} \big(t,\varphi_a^{-1}\big(x(t)\big),\Phi\big(x(t),\bm{w}(t)\big)\big) \medskip \\
x(s_1) = \varphi_a^{-1}(\bar{\bm{q}}(s_1)) \ , \ x(s_2) = \varphi_a^{-1}(\bar{\bm{q}}(s_2)) \medskip \\
\bm{c}(\bm{w}(t)) = \bm{c}\big(\varphi_a^{-1}\big(x(t)\big),\Phi\big(x(t),\bm{w}(t)\big)\big) \le 0
\end{cases}
\end{eqnarray*}
where $d\varphi_a$ is the differential of $\varphi_a$ and (recall Section \ref{localPMP})
$$
\Phi : U \times \mathbb{R}^3 \rightarrow \mathbb{R}^3 : (x,\bm{w}) \mapsto R^{\top}(x) R^{\top}_a(x) \bm{w} \ .
$$
Then, by denoting $\bar{\bm{w}}(\cdot) = R_a(\bar{x}(\cdot)) R(\bar{x}(\cdot)) \bar{\bm{u}}|_{[s_1,s_2]}(\cdot)$ with $\bar{x}(\cdot) = \varphi_a(\bar{\bm{q}}|_{[s_1,s_2]}(\cdot))$, $(\bar{x}(\cdot),\bar{\bm{w}}(\cdot))$ is optimal for (\textbf{GOGP})$_a$. \\

For what follows, recall definitions and notations in Section \ref{maxSect}. By applying the PMP to (\textbf{GOGP}) (i.e., relations \eqref{adjointSystem}-\eqref{transvCond1}), we obtain the existence of a non-positive scalar $p^0$, an absolutely continuous mapping $\bm{p} : [0,T] \rightarrow \mathbb{R}^6$ and a vector function $\bm{\mu}(\cdot) \in L^{\infty}([0,T],\mathbb{R}^3)$, with $(\bm{p}(\cdot),p^0) \neq 0$, such that, almost everywhere in $[0,T]$, the following holds
\begingroup
\begin{equation} \label{adjApp}
\displaystyle \dot{\bm{p}}(t) = -\frac{\partial H^0}{\partial \bm{q}}(t,\bar{\bm{q}}(t),\bm{p}(t),p^0,\bar{\bm{u}}(t)) - \bm{\mu}(t) \cdot \frac{\partial \bm{c}}{\partial \bm{q}}(\bar{\bm{q}}(t),\bar{\bm{u}}(t))
\end{equation}
\begin{equation} \label{maxCondApp}
\displaystyle H^0(t,\bar{\bm{q}}(t),\bm{p}(t),p^0,\bar{\bm{u}}(t)) = \underset{\bm{c}(\bar{\bm{q}}(t),\bm{u}) \le 0}{\max} H^0(t,\bar{\bm{q}}(t),\bm{p}(t),p^0,\bm{u})
\end{equation}
\begin{equation} \label{maxDerApp}
\displaystyle \frac{\partial H^0}{\partial \bm{u}}(t,\bar{\bm{q}}(t),\bm{p}(t),p^0,\bar{\bm{u}}(t)) + \bm{\mu}(t) \cdot \frac{\partial \bm{c}}{\partial \bm{u}}(\bar{\bm{q}}(t),\bar{\bm{u}}(t)) = 0
\end{equation}
\endgroup
and, in addition, conditions (\ref{slack})-(\ref{transvCond1}) hold. Since the quantity $\bm{c}\big(\bm{q},\Phi\big(\varphi_a(\bm{q}),\bm{w}\big)\big)$ does not depend on the state variable $\bm{q}$ (see \eqref{localConstraint1}, \eqref{localConstraint2}), by differentiating it w.r.t. $\bm{q}$, one obtains
$$
\frac{\partial \bm{c}}{\partial \bm{q}}(\bar{\bm{q}}(t),\bar{\bm{u}}(t)) + \frac{\partial \bm{c}}{\partial \bm{u}}(\bar{\bm{q}}(t),\bar{\bm{u}}(t)) \cdot \frac{\partial \Phi}{\partial \bm{q}}(\bar{x}(\cdot),\bar{\bm{w}}(\cdot)) = 0 \ .
$$
Moreover, by multiplying the previous expression by $\bm{\mu}(t)$ and plugging it into (\ref{maxDerApp}), we straightforwardly have
\begingroup
\small
$$
\bm{\mu}(t) \cdot \frac{\partial \bm{c}}{\partial \bm{q}}(\bar{\bm{q}}(t),\bar{\bm{u}}(t)) = \frac{\partial H^0}{\partial \bm{u}}(t,\bar{\bm{q}}(t),\bm{p}(t),p^0,\bar{\bm{u}}(t)) \cdot \frac{\partial \Phi}{\partial \bm{q}}(\bar{x}(\cdot),\bar{\bm{w}}(\cdot))
$$
\endgroup
such that, a.e. in $[s_1,s_2]$, the adjoint equation (\ref{adjApp}) becomes
\begin{multline} \label{localAdjApp}
\displaystyle \dot{\bm{p}}(t) = -\frac{\partial H^0}{\partial \bm{q}}(t,\bar{\bm{q}}(t),\bm{p}(t),p^0,\bar{\bm{u}}(t)) \\
\displaystyle -\frac{\partial H^0}{\partial \bm{u}}(t,\bar{\bm{q}}(t),\bm{p}(t),p^0,\bar{\bm{u}}(t)) \cdot \frac{\partial \Phi}{\partial \bm{q}}(\bar{x}(\cdot),\bar{\bm{w}}(\cdot)) \ .
\end{multline}
Then, by defining $p(t) = (\varphi_a^{-1})^*_{\bar{\bm{q}}(t)} \cdot \bm{p}(t)$ for every $t \in [s_1,s_2]$, it is straightforward to obtain from (\ref{localAdjApp}) and standard symplectic geometry computations (see, e.g., \cite{agrachev2013control}) that
\begingroup
\small
\begin{multline} \label{localAdjFinalApp}
\displaystyle \dot{p}(t) = -p(t) \cdot \frac{\partial}{\partial x}\Big(d\varphi_a \cdot \bm{h} \big(t,\varphi_a^{-1}(x),\Phi(x,\bm{w}))\Big)(t,\bar{x}(t),\bar{\bm{w}}(t)) \ .
\end{multline}
\endgroup
Moreover, from the properties of $\Phi$, we immediately see that the maximality condition (\ref{maxCondApp}) reads as
\begin{equation} \label{maxCondFinalApp}
\displaystyle H^0_a(t,\bar{x}(t),p(t),p^0,\bar{\bm{w}}(t)) \ge H^0_a(t,\bar{x}(t),p(t),p^0,\bm{w})
\end{equation}
for $\bm{w}$ such that $\bm{c}(\bm{w}) = \bm{c}\big(\varphi_a^{-1}\big(\bar{x}(t)\big),\Phi\big(\bar{x}(t),\bm{w}\big)\big) \le 0$ where
$$
H^0_a(t,x,p,p^0,\bm{w}) := p \cdot \Big(d\varphi_a \cdot \bm{h} \big(t,\varphi_a^{-1}(x),\Phi(x,\bm{w}))\Big)
$$
From conditions (\ref{localAdjFinalApp}), (\ref{maxCondFinalApp}), it is easily deduced that $(p(\cdot),p^0)$ is the sought multiplier for the PMP formulation related to problem (\textbf{GOGP})$_a$ (see, e.g., \cite{pontryagin1987mathematical}). Theorem \ref{mainTheo} is proved.

\subsection{Computation of Regular Controls} \label{regularApp}

In this section we compute regular controls for (\textbf{GOGP}) under Assumption \ref{assumpSimply}. We first assume that the system evolves in $(U_a,\varphi_a)$, within a non-zero measure subset $J \subseteq [0,T]$. Then, it holds $p_{\gamma}|_J(\cdot) \neq 0$ or $p_{\chi}|_J(\cdot) \neq 0$ (see Section \ref{controlSection}). \\

If $p^a_v|_J(\cdot) = 0$, by definition $C_a|_J(\cdot) = D_a|_J(\cdot) = 0$ and then, from (\ref{firstHam}) and the Cauchy-Schwarz inequality, we obtain
\begin{equation*}
\displaystyle w_2 = \frac{\sin(\alpha_{\max}) p_{\gamma}}{\sqrt{p^2_{\gamma} + \frac{p^2_{\chi}}{\cos^2(\gamma)}}} \quad , \quad w_3 = \frac{\sin(\alpha_{\max}) p_{\chi}}{\cos(\gamma) \sqrt{p^2_{\gamma} + \frac{p^2_{\chi}}{\cos^2(\gamma)}}} \ .
\end{equation*}
Therefore, $w_1 = \sqrt{1 - (w^2_2 + w^2_3)}$ thanks to constraint $c_1$. \\

We analyze now the harder case $p^a_v|_J(\cdot) \neq 0$. Denote $\lambda = p_{\gamma} \omega$, $\rho = p_{\chi} \frac{\omega}{\cos(\gamma)}$. In the following, we apply the Karush-Kuhn-Tucker conditions to (\ref{firstHam}). For this, we first remark that any optimum for (\ref{firstHam}) satisfies $w_1 > 0$. Moreover, if the constraints in (\ref{firstHam}) were active at the optimum, then this point would satisfy $\bm{w} \in S^2$, $w^2_2 + w^2_3 = \sin^2(\alpha_{\max})$, and consequently, the gradients of these constraints evaluated at the optimum would satisfy the linear independence constraint qualification (see, e.g., \cite{NoceWrig06}). By applying the Karush-Kuhn-Tucker conditions to (\ref{firstHam}) (without considering $w_1 \ge 0$, thanks to what we said above), we infer the existence of a non-zero multiplier $(\eta_1,\eta_2) \in \mathbb{R} \times \mathbb{R}_+$ which satisfies
\begingroup
\small
\begin{eqnarray*}
\begin{cases}
C_a - 2 \eta_1 w_1 = 0 \quad , \quad 2 (\eta_1 + \eta_2 + D_a) w_2 - \lambda = 0 \medskip \\
2 (\eta_1 + \eta_2 + D_a) w_3 - \rho = 0 \ , \ \eta_2 (w^2_2 + w^2_3 - \sin^2(\alpha_{\max})) = 0 \ .
\end{cases}
\end{eqnarray*}
\endgroup
Since either $\lambda \neq 0$ or $\rho \neq 0$, one necessarily has $\eta_1 + \eta_2 + D_a \neq 0$ so that the optimal control satisfies $\rho w_2 = \lambda w_3$. We proceed considering $\lambda \neq 0$, i.e., $w_3 = (\rho/\lambda) w_2$. The problem is reduced to the study of the following optimization
\begingroup
\begin{multline*}
\max \Bigg\{ C_a w_1 - \bigg(1 + \frac{\rho^2}{\lambda^2}\bigg) (D_a w^2_2 - \lambda w_2) \ \mid \\
w^2_1 + \bigg(1 + \frac{\rho^2}{\lambda^2}\bigg) w^2_2 = 1 \ , \ \bigg(1 + \frac{\rho^2}{\lambda^2}\bigg) w^2_2 \le \sin^2(\alpha_{\max}) \Bigg\} \ .
\end{multline*}
\endgroup
In other words, we seek points $(w_1,w_2)$ such that the relations
\begingroup
\begin{multline} \label{parabola}
w_1 = \frac{1}{C_a} \bigg(1 + \frac{\rho^2}{\lambda^2}\bigg) (D_a w^2_2 - \lambda w_2) + \frac{C}{C_a} \ , \\
\ w^2_1 + \bigg(1 + \frac{\rho^2}{\lambda^2}\bigg) w^2_2 = 1 \ , \ \bigg(1 + \frac{\rho^2}{\lambda^2}\bigg) w^2_2 \le \sin^2(\alpha_{\max})
\end{multline}
\endgroup
hold with the largest possible $C \in \mathbb{R}$. Several cases occur:
\begin{itemize}
\item $C_a = 0 \ $:

Since $D_a \neq 0$, this case results in the maximization of a parabola under box constraints. By denoting $A = - \Big(1 + \frac{\rho^2}{\lambda^2}\Big) D_a$, $B = \Big(1 + \frac{\rho^2}{\lambda^2}\Big) \lambda$ and $D = \frac{|\lambda| \sin(\alpha_{\max})}{\sqrt{\lambda^2 + \rho^2}}$, we maximize $A w^2_2 + B w_2$ such that $-D \le w_2 \le D$. Then, one has $w_1 = \sqrt{1 - \Big(1 + \frac{\rho^2}{\lambda^2}\Big) w^2_2}$, where:
\begin{itemize}
\item $w_2 = -D$ if $A > 0$, $B < 0$ or $A > 0$, $B < -2 |A| D$;
\item $w_2 = -\frac{B}{2 A}$ if $A > 0$, $-2 |A| D \le B \le 2 |A| D$;
\item $w_2 = D$ if $A > 0$, $B > 0$ or $A < 0$, $B > 2 |A| D$.
\end{itemize}
\item $C_a > 0 \ $:

The optimum is given by the contact point between the parabola and the ellipse given in (\ref{parabola}) that lies in the positive half-plane $w_1 > 0$. Under Assumption \ref{assumpSimply}, this is given by matching the first derivatives of these curves. More specifically, this provides $w_1 = \sqrt{1 - \frac{\lambda^2 + \rho^2}{(C_a + 2 D_a)^2}}$, $w_2 = \frac{\lambda}{C_a + 2 D_a}$ if $\frac{\lambda^2 + \rho^2}{(C_a + 2 D_a)^2} \le \sin^2(\alpha_{\max})$. However, saturations may arise, i.e., $w_1 = \cos(\alpha_{\max})$ and $w_2 = -\frac{|\lambda| \sin(\alpha_{\max})}{\displaystyle \sqrt{\lambda^2 + \rho^2}}$ if $\frac{\lambda}{C_a + 2 D_a} < -\frac{|\lambda| \sin(\alpha_{\max})}{\sqrt{\lambda^2 + \rho^2}}$, or $w_2 = \frac{|\lambda| \sin(\alpha_{\max})}{\sqrt{\lambda^2 + \rho^2}}$ if $\frac{\lambda}{C_a + 2 D_a} > \frac{|\lambda| \sin(\alpha_{\max})}{\sqrt{\lambda^2 + \rho^2}}$.
\item $C_a < 0 \ $:

In this case, since $w_1 > 0$, the optimum becomes the point of intersection beetwen the parabola and the upper part of the ellipse given in (\ref{parabola}) for which $C$ takes the maximum value. Only saturations are allowed. Indeed, by studying the position of the minimum of the parabola, we obtain $w_1 = \cos(\alpha_{\max})$ and $w_2 = -\frac{|\lambda| \sin(\alpha_{\max})}{\sqrt{\lambda^2 + \rho^2}}$ if $\frac{\lambda}{D_a} > 0$, or $w_2 = \frac{|\lambda| \sin(\alpha_{\max})}{\sqrt{\lambda^2 + \rho^2}}$ if $\frac{\lambda}{D_a} < 0$.
\end{itemize}
Clearly, a similar procedure holds when $\rho \neq 0$, $w_2 = (\lambda/\rho) w_3$. \\

At this step, we have found the optimal strategy in the regular case for the first local chart $(U_a,\varphi_a)$. By the similarity of (\ref{firstHam}) and (\ref{secondHam}), similar results hold true for the local control $\bm{z}$ using instead the second local chart $(U_b,\varphi_b)$, for which $\lambda$ and $\rho$ are replaced respectively by $p_{\theta} \omega$ and by $-p_{\phi} \frac{\omega}{\cos(\theta)}$. \\

We have found the behavior of any regular controls.

\subsection{Computation of Nonregular Controls} \label{singularApp}

In this section we compute nonregular optimal controls for (\textbf{GOGP}), under Assumption \ref{assM} and Assumption \ref{assVel}, within a non-zero measure subset $J \subseteq [0,T]$. In what follows, we will need the adjoint equations related to (\textbf{GOGP})$_a$. These come from applying the PMP for problems with pure control constraints (see, e.g., \cite{pontryagin1987mathematical}) to (\textbf{GOGP})$_a$ and are listed below:
\begingroup
\footnotesize
\begin{multline*}
\displaystyle \dot{p}^a_r = p^a_L \frac{v}{r^2} \cos(\gamma) \cos(\chi) + p^a_l \frac{v}{r^2} \frac{\cos(\gamma) \sin(\chi)}{\cos(L)} + p_{\gamma} \bigg( \frac{v c_m}{h_r} w_2 \\
+ \frac{v}{r^2} \cos(\gamma) + \frac{\partial g}{\partial r} \frac{\cos(\gamma)}{v} \bigg) + p_{\chi} \bigg( \frac{v c_m}{h_r \cos(\gamma)} w_3 + \frac{v}{r^2} \cos(\gamma) \sin(\chi) \tan(L) \bigg) \medskip \\
+ p^a_v \bigg( \frac{\partial g}{\partial r} \sin(\gamma) - \frac{v^2}{h_r} \big( d + \eta c_m (w^2_2 + w^2_3) \big) \bigg)
\end{multline*}
$$
\displaystyle \dot{p^a}_L = -p^a_l \frac{v}{r} \frac{\cos(\gamma) \sin(\chi) \tan(L)}{\cos(L)} - p_{\chi} \frac{v}{r} \frac{\cos(\gamma) \sin(\chi)}{\cos^2(L)} \ , \ \dot{p}^a_l = 0
$$
\begin{multline*}
\displaystyle \dot{p}^a_v = -p^a_r \sin(\gamma) - p^a_L \frac{\cos(\gamma) \cos(\chi)}{r} - p^a_l \frac{\cos(\gamma) \sin(\chi)}{r \cos(L)} \\
\displaystyle + 2 p^a_v v \big( d + \eta c_m (w^2_2 + w^2_3) \big) + p_{\gamma} \bigg( \frac{\omega}{v} w_2 - \frac{\cos(\gamma)}{r} - \frac{g}{v^2} \cos(\gamma) \bigg) \\
\displaystyle  + p_{\chi} \bigg( \frac{\omega}{v} \frac{w_3}{\cos(\gamma)} - \frac{\cos(\gamma) \sin(\chi) \tan(L)}{r} \bigg)
\end{multline*}
\begin{multline*}
\displaystyle \dot{p}_{\gamma} = -p^a_r v \cos(\gamma) + p^a_L \frac{v}{r} \sin(\gamma) \cos(\chi) + p^a_l \frac{v}{r} \frac{\sin(\gamma) \sin(\chi)}{\cos(L)} + p_{\gamma} \bigg( \frac{v}{r} \\
\displaystyle - \frac{g}{v} \bigg) \sin(\gamma) + p^a_v g \cos(\gamma) + p_{\chi} \bigg( \frac{v}{r} \sin(\gamma) \sin(\chi) \tan(L) - \frac{\omega \sin(\gamma)}{\cos^2(\gamma)} w_3 \bigg)
\end{multline*}
$$
\displaystyle \dot{p}_{\chi} = p^a_L \frac{v}{r} \cos(\gamma) \sin(\chi) - p^a_l \frac{v}{r} \frac{\cos(\gamma) \cos(\chi)}{\cos(L)} - p_{\chi} \frac{v}{r} \cos(\gamma) \cos(\chi) \tan(L) .
$$
\endgroup

The important result that allows us to work out explicit expressions for nonregular controls consists of showing that, under Assumption \ref{assM}, it holds $p^a_v|_J(\cdot) \neq 0$, $p^b_v|_J(\cdot) \neq 0$, i.e., nonregular controls are not degenerate. This arises as follows.

\begin{lemma} \label{propSingCont}
Suppose $p_{\gamma}|_J(\cdot) = p_{\chi}|_J(\cdot) = 0$ (as well as $p_{\theta}|_J(\cdot) = p_{\phi}|_J(\cdot) = 0$), i.e., nonregular controls appear. Then, under Assumption \ref{assM}, $p^a_v|_J(\cdot) \neq 0$ (as well as $p^b_v|_J(\cdot) \neq 0$).
\end{lemma}

\textit{Proof:} We prove the statement considering the first local chart $(U_a,\varphi_a)$. For the second local chart, similar computations hold. By contradiction, suppose that $p_{\gamma}|_J(\cdot) = p_{\chi}|_J(\cdot) = p^a_v|_J(\cdot) = 0$. From the adjoint equations of coordinates $p^a_v$, $p_{\gamma}$ and $p_{\chi}$ (given above) restricted to $J$, we obtain
\begingroup
\footnotesize
$$
\left( \arraycolsep=1.5pt \begin{array}{ccc}
-v \cos(\gamma) & \displaystyle \frac{v}{r} \sin(\gamma) \cos(\chi) & \displaystyle \frac{v}{r} \frac{\sin(\gamma) \sin(\chi)}{\cos(L)} \medskip  \\
0 & \displaystyle \frac{v}{r} \cos(\gamma) \sin(\chi) & \displaystyle -\frac{v}{r} \frac{\cos(\gamma) \cos(\chi)}{\cos(L)} \medskip \\
-\sin(\gamma) & \displaystyle \frac{\cos(\gamma) \cos(\chi)}{r} & \displaystyle \frac{\cos(\gamma) \sin(\chi)}{r \cos(L)} \medskip  \\
\end{array} \right) \left( \begin{array}{c}
p^a_r \medskip \\
p^a_L \medskip \\
p^a_l
\end{array} \right) =
\left( \begin{array}{c}
0 \medskip \\
0 \medskip \\
0
\end{array} \right) \ .
$$
\endgroup
The determinant of the matrix is $\frac{v^2 \cos(\gamma)}{r^2 \cos(L)} \neq 0$ and then $(p^a_r,p^a_L,p^a_l)|_J(\cdot) = 0$. This implies that the adjoint vector is zero everywhere in $[0,T]$. Assumption \ref{assM}, the transversality conditions and $\bm{p}(\cdot) \equiv 0$ (from Theorem \ref{mainTheo}) give $p^0 = 0$, raising thus a contradiction because it holds $(\bm{p}(\cdot),p^0) \neq 0$. $_{\Box}$ \\

\subsubsection{First Local Chart Representation}

We start assuming that the system evolves in the first local chart $(U_a,\varphi_a)$, within a non-zero mesure subset $J \subseteq [0,T]$. Our objective consists in studying (\ref{singularHam}). Thanks to Lemma \ref{propSingCont}, from now on,  we assume $p_{\gamma}|_J(\cdot) = p_{\chi}|_J(\cdot) = 0$, $p^a_v|_J(\cdot) \neq 0$ and, when clear from the context, we skip the dependence on $t$ to keep better readability. Moreover, we introduce the following local forms for the dynamics of (\textbf{GOGP})$_a$ (recall Section \ref{firstChart}):
\begingroup
\footnotesize
\begin{multline*}
\displaystyle X(t,\bm{r},\bm{v}) := v\sin(\gamma) \frac{\partial}{\partial r} + \frac{v}{r} \cos(\gamma) \cos(\chi) \frac{\partial}{\partial L} + \frac{v}{r} \frac{\cos(\gamma) \sin(\chi)}{\cos(L)} \frac{\partial}{\partial l} \\
\displaystyle - \left(d v^2 + g \sin(\gamma) \right) \frac{\partial}{\partial v} + \left(\frac{v}{r} - \frac{g}{v}\right) \cos(\gamma) \frac{\partial}{\partial \gamma} + \frac{v}{r} \cos(\gamma) \sin(\chi) \tan(L) \frac{\partial}{\partial \chi}
\end{multline*}
\begin{equation*}
\displaystyle Y_1(t,\bm{r},\bm{v}) := \frac{f_T}{m} \frac{\partial}{\partial v} \quad , \quad Y_Q(t,\bm{r},\bm{v}) := -\eta c_m v^2 \frac{\partial}{\partial v}
\end{equation*}
\begin{equation*}
\displaystyle Y_2(t,\bm{r},\bm{v}) := \omega \frac{\partial}{\partial \gamma} \quad , \quad Y_3(t,\bm{r},\bm{v}) := \frac{\omega}{\cos(\gamma)} \frac{\partial}{\partial \chi} \ .
\end{equation*}
\endgroup
The Lie bracket of two vector fields $X$, $Y$ is defined as the derivation $[X,Y](f) := X(Yf) - Y(Xf)$, $f \in C^{\infty}$ (see, e.g., \cite{agrachev2013control}). The following classical result holds (see, e.g., \cite{bonnard2003optimal}).
\begin{lemma} \label{lemmaLie}
Using the first local chart $(U_a,\varphi_a)$, for times $t \in J$ such that $(\bm{r},\bm{v})(t)$ lies within $U_a$, the following holds:
\begingroup
\small
\begin{multline} \label{lie1}
\displaystyle \frac{d}{d t} \big\langle \bm{p} , Y_2 \big\rangle = \big\langle \bm{p} , \frac{\partial}{\partial t} Y_2 \big\rangle + \big\langle \bm{p} , [X,Y_2] \big\rangle + w_1 \big\langle \bm{p} , [Y_1,Y_2] \big\rangle \\
\displaystyle + w_3 \big\langle \bm{p} , [Y_3,Y_2] \big\rangle + (w^2_2 + w^2_3) \big\langle \bm{p} , [Y_Q,Y_2] \big\rangle
\end{multline}
\begin{multline} \label{lie2}
\displaystyle \frac{d}{d t} \big\langle \bm{p} , Y_3 \big\rangle = \big\langle \bm{p} , \frac{\partial}{\partial t} Y_3 \big\rangle + \big\langle \bm{p} , [X,Y_3] \big\rangle + w_1 \big\langle \bm{p} , [Y_1,Y_3] \big\rangle \\
\displaystyle + w_2 \big\langle \bm{p} , [Y_2,Y_3] \big\rangle + (w^2_2 + w^2_3) \big\langle \bm{p} , [Y_Q,Y_3] \big\rangle
\end{multline}
\begin{multline} \label{lie3}
\displaystyle \frac{d}{d t} \big\langle \bm{p} , [X,Y_2] \big\rangle = \big\langle \bm{p} , \frac{\partial}{\partial t} [X,Y_2] \big\rangle + \big\langle \bm{p} , [X,[X,Y_2]] \big\rangle \\
\displaystyle + w_1 \big\langle \bm{p} , [Y_1,[X,Y_2]] \big\rangle + w_2 \big\langle \bm{p} , [Y_2,[X,Y_2]] \big\rangle \\
\displaystyle + w_3 \big\langle \bm{p} , [Y_3,[X,Y_2]] \big\rangle + (w^2_2 + w^2_3) \big\langle \bm{p} , [Y_Q,[X,Y_2]] \big\rangle
\end{multline}
\begin{multline} \label{lie4}
\displaystyle \frac{d}{d t} \big\langle \bm{p} , [X,Y_3] \big\rangle = \big\langle \bm{p} , \frac{\partial}{\partial t} [X,Y_3] \big\rangle + \big\langle \bm{p} , [X,[X,Y_3]] \big\rangle \\
\displaystyle + w_1 \big\langle \bm{p} , [Y_1,[X,Y_3]] \big\rangle + w_2 \big\langle \bm{p} , [Y_2,[X,Y_3]] \big\rangle \\
\displaystyle + w_3 \big\langle \bm{p} , [Y_3,[X,Y_3]] \big\rangle + (w^2_2 + w^2_3) \big\langle \bm{p} , [Y_Q,[X,Y_3]] \big\rangle \ .
\end{multline}
\endgroup
\end{lemma}

The idea that we develop here seeks explicit expressions for the optimal controls $\bm{w}(\cdot)$ by analyzing expressions (\ref{lie1})-(\ref{lie4}). Our strategy is based on the following remarks, which come from symbolic Lie bracket computations on the local fields:
\begin{enumerate}
\item[(A)] $[Y_1,Y_2]$, $[Y_Q,Y_2]$ are proportional to $\frac{\partial}{\partial \gamma}$;
\item[(B)] $[Y_1,Y_3]$, $[Y_2,Y_3]$, $[Y_Q,Y_3]$, $[Y_2,[X,Y_3]]$ lie along $\frac{\partial}{\partial \chi}$;
\item[(C)] When $p_{\gamma}|_J(\cdot) = p_{\chi}|_J(\cdot) = 0$, then $\big\langle \bm{p} , [X,[X,Y_3]] \big\rangle$, $\big\langle \bm{p} , [Y_1,[X,Y_3]] \big\rangle$, $\big\langle \bm{p} , [Y_Q,[X,Y_3]] \big\rangle$ lie along $\dot{p}_{\chi}$;
\item[(D)] When $p_{\gamma}|_J(\cdot) = p_{\chi}|_J(\cdot) = 0$, $\big\langle \bm{p} , \frac{\partial}{\partial t} [X,Y_2] \big\rangle$ lies along $\big\langle \bm{p} , [X,Y_2] \big\rangle$ and $\big\langle \bm{p} , \frac{\partial}{\partial t} [X,Y_3] \big\rangle$ lies along $\big\langle \bm{p} , [X,Y_3] \big\rangle$.
\end{enumerate}
From $p_{\gamma}|_J(\cdot) = p_{\chi}|_J(\cdot) = 0$, (A) and (B) applied to (\ref{lie1}) and (\ref{lie2}) give $\big\langle \bm{p} , [X,Y_2] \big\rangle \big|_J = \big\langle \bm{p} , [X,Y_3] \big\rangle \big|_J = 0$. These expressions, plugged into (\ref{lie4}) using (B), (C) and (D), lead to
\begin{equation} \label{exprSing}
w_3 \big\langle \bm{p} , [Y_3,[X,Y_3]] \big\rangle = 0 \ , \ \textnormal{in } J \ .
\end{equation}
Seeking explicit expressions for the nonregular controls from (\ref{exprSing}) becomes a hard and tedious task when $\big\langle \bm{p} , [Y_3,[X,Y_3]] \big\rangle = 0$. This because more many time derivatives are required, which provide complex expressions of Lie brackets. In this situation, the environmental conditions concerning the feasibility of (\textbf{GOGP}) (represented by Assumption \ref{assVel}) play an important role in making these further time derivatives of Lie brackets not necessary for our purpose. Indeed, we have the following:

\begin{lemma} \label{lemmaDuality}
Assume that Assumption \ref{assVel} holds. Then, one has $\big\langle \bm{p} , [Y_3,[X,Y_3]] \big\rangle \neq 0$ almost everywhere in $J$.
\end{lemma}
\textit{Proof:} By contradiction, suppose that $\big\langle \bm{p} , [Y_3,[X,Y_3]] \big\rangle = 0$ a.e. within $J$. This implies that $\cos(\chi) p^a_L + \frac{\sin(\chi)}{\cos(L)} p^a_l = 0$ a.e. within $J$. The previous expression, combined with the adjoint equation for $p_{\chi}$ (given above), gives $p^a_L|_J(\cdot) = p^a_l|_J(\cdot) = 0$. On the other hand, from the adjoint equation of $p_{\gamma}$ (see above), we have $(v p^a_r - g p^a_v)|_J(\cdot) = 0$. Combining this expression with its derivative w.r.t. time and imposing $p^a_v|_J(\cdot) \neq 0$ lead to
$$
v^4 + 3 g(\bm{r}) h_r v^2 - g(\bm{r}) h_r \left( \frac{f_T w_1}{m (d + \eta c_m (w^2_2 + w^2_3))} \right) = 0 \ .
$$
First of all, if $f_T = 0$ a contradiction arises immediately. The only physically meaningful solution for this equation is
\begingroup
\small
$$
v = \sqrt{\frac{3}{2} g(\bm{r}) h_r} \sqrt{ \sqrt{1 + \frac{4}{9} \frac{1}{g(\bm{r}) h_r} \left( \frac{f_T w_1}{m (d + \eta c_m (w^2_2 + w^2_3))} \right) } - 1 }
$$
\endgroup
and, from $w_1 \in [0,1]$, Assumption \ref{assVel} gives a contradiction. $_{\Box}$

The previous results allow us to reformulate (\ref{singularHam}) as
\begingroup
\small
\begin{equation*}
\displaystyle (w_1,w_2) = \argmax \Big\{ C_a w_1 - D_a w^2_2 \mid w^2_1 + w^2_2 = 1, w^2_2 \le \sin^2(\alpha_{\max}) \Big\}
\end{equation*}
\endgroup
that we can solve. Remark that $D_a \neq 0$, $C_a \neq 0$ iff $f_T \neq 0$.

Suppose first that $C_a = 0$. In this case, it is clear that the maximization problem above is solved by $w_1 = 1$, $w_2 = 0$ if $D_a > 0$ and $w_1 = \cos(\alpha_{\max})$, $w^2_2 = \sin^2(\alpha_{\max})$ if $D_a < 0$. Let now $C_a \neq 0$. Exploiting a graphical study, it is not difficult to see that the solutions are now given by $w_1 = 1$, $w_2 = 0$ if $C_a > 0$ and $w_1 = \cos(\alpha_{\max})$, $w^2_2 = \sin^2(\alpha_{\max})$ if $C_a < 0$.

To conclude, it remains to establish the value of the coordinate $w_2$ when $w_1 = \cos(\alpha_{\max})$ and $w^2_2 = \sin^2(\alpha_{\max})$. For this, we may use expression (\ref{lie3}). Indeed, it is clear that, when $\big\langle \bm{p} , [Y_2,[X,Y_2]] \big\rangle \neq 0$, it holds (recall statements (A)-(D))
\begingroup
\begin{multline*}
w_2 \displaystyle = -\frac{\big\langle \bm{p} , [X,[X,Y_2]] \big\rangle}{\big\langle \bm{p} , [Y_2,[X,Y_2]] \big\rangle} - w_1 \frac{\big\langle \bm{p} , [Y_1,[X,Y_2]] \big\rangle}{\big\langle \bm{p} , [Y_2,[X,Y_2]] \big\rangle} \\
- w^2_2 \frac{\big\langle \bm{p} , [Y_Q,[X,Y_2]] \big\rangle}{\big\langle \bm{p} , [Y_2,[X,Y_2]] \big\rangle} \ .
\end{multline*}
\endgroup
If instead $\big\langle \bm{p} , [Y_2,[X,Y_2]] \big\rangle = 0$ a.e. in $J$, then, suppose that $\big\langle \bm{p} , [Y_2,[Y_2,[X,Y_2]]] \big\rangle \neq 0$. Differentiating as done in (\ref{lie3}), (\ref{lie4}), by using the same arguments as above we have
\begingroup
\begin{multline*}
w_2 \displaystyle = -\frac{\Big\langle \bm{p} , [Y_2,[X,[X,Y_2]]] \Big\rangle}{\Big\langle \bm{p} , [Y_2,[Y_2,[X,Y_2]]] \Big\rangle} - w_1 \frac{\Big\langle \bm{p} , [Y_2,[Y_1,[X,Y_2]]] \Big\rangle}{\Big\langle \bm{p} , [Y_2,[Y_2,[X,Y_2]]] \Big\rangle} \\
- w^2_2 \frac{\Big\langle \bm{p} , [Y_2,[Y_Q,[X,Y_2]]] \Big\rangle}{\Big\langle \bm{p} , [Y_2,[Y_2,[X,Y_2]]] \Big\rangle} \ .
\end{multline*}
\endgroup
We can prove that actually one between the two previous formulas always holds, giving then the sought conclusion.
\begin{lemma}
Under Assumption \ref{assVel}, almost everywhere in $J$:
$$
\big\langle \bm{p} , [Y_2,[X,Y_2]] \big\rangle \neq 0 \quad \textnormal{or} \quad \big\langle \bm{p} , [Y_2,[Y_2,[X,Y_2]]] \big\rangle \neq 0 \ .
$$
\end{lemma}
\textit{Proof:} By contradiction, suppose that $\big\langle \bm{p} , [Y_2,[X,Y_2]] \big\rangle = 0$ and $\big\langle \bm{p} , [Y_2,[Y_2,[X,Y_2]]] \big\rangle = 0$ a.e. in $J$. From these, one recovers respectively the following two expressions
\begingroup
\footnotesize
$$
\bigg( \sin(\gamma) p^a_r + \frac{\cos(\gamma) \cos(\chi)}{r} p^a_L + \frac{\cos(\gamma) \sin(\chi)}{r \cos(L)} p^a_l - \frac{g \sin(\gamma)}{v} p^a_v \bigg) \bigg|_J(\cdot) = 0
$$
$$
\bigg( \cos(\gamma) p^a_r - \frac{\sin(\gamma) \cos(\chi)}{r} p^a_L - \frac{\sin(\gamma) \sin(\chi)}{r \cos(L)} p^a_l - \frac{g \cos(\gamma)}{v} p^a_v \bigg) \bigg|_J(\cdot) = 0
$$
\endgroup
which lead to $\cos(\chi) p^a_L + \frac{\sin(\chi)}{\cos(L)} p^a_l = 0$ a.e. within $J$. This relation, combined with the adjoint equation of $p_{\chi}$ (see above), gives $p^a_L|_J(\cdot) = p^a_l|_J(\cdot) = 0$. On the other hand, the adjoint equation of $p_{\gamma}$ (above) provides $(v p^a_r - g p^a_v)|_J(\cdot) = 0$. Then, as in the proof of Lemma \ref{lemmaDuality}, a contradiction arises. $_{\Box}$ \\

\subsubsection{Second Local Chart Representation}

The approach proposed in the previous section is no more exploitable when using the second local chart $(U_b,\varphi_b)$ and problem (\ref{singularHam2}). Indeed, the terms of the gravity and the curvature of the Earth contained in (\ref{dynSecond}) make the computations on the Lie algebra generated by the local fields hard to treat. However, we can still recover nonregular arcs by proceeding as follows.

Thanks to the previous computation, we know the explicit expressions for nonregular controls for every point in $U_a$. Then, it is enough to compute possible nonregular controls for trajectories in $U_b \setminus U_a$. From (\ref{frame1}), (\ref{frame2}), one sees that these trajectories lie in the following submanifold of $\mathbb{R}^6 \setminus \{ 0 \}$
$$
S_b := \left\{ (\bm{r},\bm{v}) \in \mathbb{R}^6 \setminus \{ 0 \} \ \mid \ \bm{v} \ \mathbin{\!/\mkern-5mu/\!} \ \bm{r} \right\}
$$
which corresponds, by forcing the coordinates of the chart $(U_b,\varphi_b)$, to points such that $\theta = 0$, $\phi = 0$ or $\theta = 0$, $\phi = \pi$. Following the previous argument, suppose that there exists a non-zero measure subset $J \subseteq [0,T]$ for which the optimal trajectory $(\bm{r},\bm{v})(\cdot)$ arisen from a nonregular control $\bm{u}(\cdot)$ is such that $(\bm{r},\bm{v})(t) \in S_b$ for every $t \in J$. In particular, suppose that $\theta|_J(\cdot) = 0 \ , \ \phi|_J(\cdot) = 0$ or $\phi|_J(\cdot) = \pi$. Then, almost everywhere in $J$, the trajectory $(\bm{r},\bm{v})(\cdot)$ satisfies
\begingroup
\begin{eqnarray*}
\begin{cases}
\dot{r} = \pm v \ , \ \dot{L} = 0 \ , \ \dot{l} = 0  \ , \ \dot{\theta} = \displaystyle \omega z_2 \ , \ \dot{\phi} = \displaystyle -\omega z_3 \medskip \\
\dot{v} = \displaystyle \frac{f_T}{m} z_1 -\left(d + \eta c_m (z^2_2 + z^2_3) \right) v^2 \pm g \ .
\end{cases}
\end{eqnarray*}
\endgroup
Since the values of $\theta$ and $\phi$ remain the same along $J$, their derivative w.r.t. the time must be zero. Therefore, almost everywhere in $J$, any nonregular control satisfies $z_1|_J(\cdot) = 1$, $z_2|_J(\cdot) = 0$ and $z_3|_J(\cdot) = 0$. This concludes our analysis.

\begin{IEEEbiography}[{\includegraphics[width=1in,height=1.25in,clip,keepaspectratio]{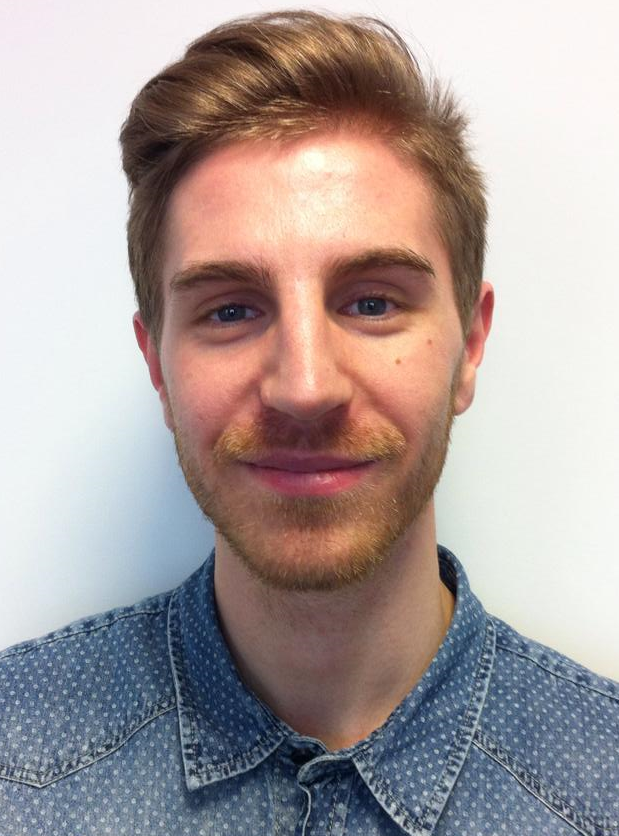}}]{Riccardo Bonalli} obtained his MSc in Mathematical Engineering from Politecnico di Milano, Italy, in 2014, and his PhD in applied mathematics from Sorbonne Universit\'e, France, in 2018, in collaboration with ONERA - The French Aerospace Lab, France. He is now postdoctoral researcher at the Department of Aeronautics and Astronautics, at Stanford University, California. His main research interests concern the theoretical and numerical optimal control with applications in aerospace engineering and robotics.
\end{IEEEbiography}
\begin{IEEEbiography}[{\includegraphics[width=1in,height=1.25in,clip,keepaspectratio]{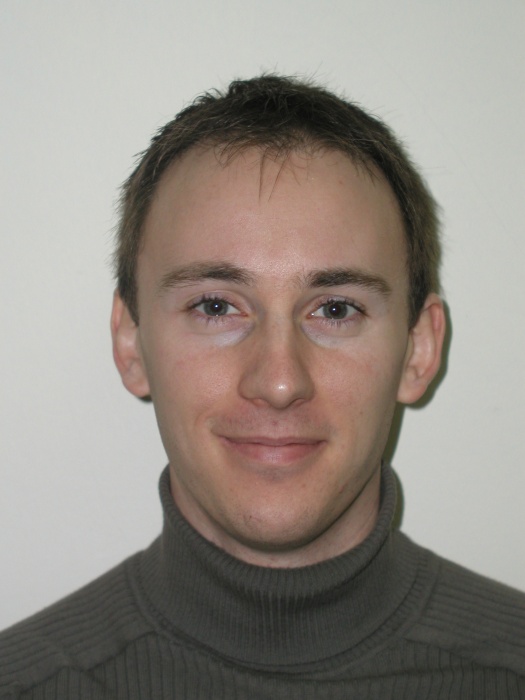}}]{Bruno H\'{e}riss\'{e}} received the Engineering degree and the Master degree from the \'Ecole Sup\'erieure d'\'Electricit\'e (SUPELEC), Paris, France, in 2007. After three years of research with CEA List, he received the Ph.D. degree in robotics from the University of Nice Sophia Antipolis, Sophia Antipolis, France, in 2010. Since 2011, he has been a Research Engineer with ONERA, the French Aerospace Lab, Palaiseau, France. His research interests include optimal control and vision-based control with applications in aerospace sytems and aerial robotics.
\end{IEEEbiography}
\begin{IEEEbiography}[{\includegraphics[width=1in,height=1.25in,clip,keepaspectratio]{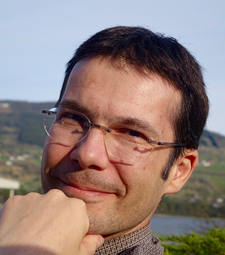}}]{Emmanuel Tr\'{e}lat} was born in 1974. He is currently full professor at Sorbonne Universit\'e (Paris 6). He is the director of the Fondation Sciences Math\'{e}matiques de Paris. He is editor in chief of the journal ESAIM: Control, Optimization and Calculus of Variations, and is associated editor of many other journals. He has been awarded the SIAM Outstanding Paper Prize (2006), Maurice Audin Prize (2010), Felix Klein Prize (European Math. Society, 2012), Blaise Pascal Prize (french Academy of Science, 2014), Big Prize Victor Noury (french Academy of Science, 2016). His research interests range over control theory in finite and infinite dimension, optimal control, stabilization, geometry, numerical analysis, with a special interest to optimal control applied to aerospace.
\end{IEEEbiography}


\end{document}